\documentclass[12pt,fleqn]{article}
\usepackage{amssymb,latexsym,amsmath,amsfonts}
\usepackage{graphicx}
\usepackage{comment}

\usepackage[thinlines]{easymat}
\usepackage{overpic}

\usepackage{a4wide}{\large }
\usepackage{rotating}

    \numberwithin{equation}{section}

    \def\ds{\displaystyle}
    
    \def\Re{{\rm Re \,  }}
    \def\Im{{\rm Im \, }}
    \def\Ai{{\rm Ai  }}

    \DeclareMathOperator*{\Tr}{Tr}

    \newtheorem{theorem}{Theorem}[section]
    \newtheorem{lemma}[theorem]{Lemma}

    \newtheorem{Definition}[theorem]{Definition}
    
    \newtheorem{Remark}[theorem]{Remark}
    \newenvironment{remark}{\begin{Remark}\rm}{\end{Remark}}
    \newtheorem{Example}[theorem]{Example}

    \newenvironment{proof}%
    {\rm \trivlist \item[\hskip \labelsep{\bf Proof. }]}%
    {\hspace*{\fill}$\Box$\endtrivlist}
    {\rm \trivlist \item[\hskip \labelsep{\bf Proof}]}%
    {\hspace*{\fill}$\Box$\endtrivlist}

\begin{document}

\title{Asymptotics for a special solution of \\
the thirty fourth Painlev\'e equation}
\author{A.R. Its \\
{\em Department of Mathematical Sciences} \\
{\em Indiana University -- Purdue University Indianapolis} \\
{\em Indianapolis IN 46202-3216, U.S.A.} \\
itsa@math.iupui.edu
 \\[15pt]
A.B.J. Kuijlaars  \\
{\em Department of Mathematics} \\
{\em Katholieke Universiteit Leuven} \\
{\em Celestijnenlaan 200B} \\
{\em 3001 Leuven, Belgium} \\
arno.kuijlaars@wis.kuleuven.be
\\[15pt]
and
\\[15pt]
J. \"{O}stensson  \\
{\em Department of Mathematics} \\
{\em Uppsala University} \\
{\em 751 06 Uppsala, Sweden} \\
ostensson@math.uu.se}

\date{November 24, 2008}

\maketitle

\newpage
\begin{abstract}
In a previous paper we studied the double scaling limit of unitary
random matrix ensembles of the form $Z_{n,N}^{-1} |\det
M|^{2\alpha} e^{-N \Tr V(M)} dM$ with $\alpha
> -1/2$. The factor $|\det M|^{2\alpha}$ induces critical
eigenvalue behavior near the origin. Under the assumption that the
limiting mean eigenvalue density associated with $V$ is regular,
and that the origin is a right endpoint of its support, we
computed the limiting eigenvalue correlation kernel in the double
scaling limit as $n, N \to \infty$ such that $n^{2/3}(n/N-1) =
O(1)$ by using the Deift-Zhou steepest descent method for the
Riemann-Hilbert problem for polynomials on the line orthogonal
with respect to the weight $|x|^{2\alpha} e^{-NV(x)}$. Our main
attention was on the construction of a local parametrix near the
origin by means of the $\psi$-functions associated with a
distinguished solution $u_{\alpha}$ of the Painlev\'e XXXIV equation.
This solution is related to a particular solution of the
Painlev\'e II equation, which however is different from the usual
Hastings-McLeod solution. In this paper we compute the asymptotic
behavior of $u_{\alpha}(s)$ as $s \to \pm \infty$. We conjecture
that this asymptotics characterizes $u_{\alpha}$ and we present supporting
arguments based on the asymptotic analysis of a one-parameter family
of solutions of the Painlev\'e XXXIV equation
which includes $u_{\alpha}$. We identify this family
as the family of {\it tronqu\'ee} solutions of the thirty fourth Painlev\'e equation.
\end{abstract}

\section{Introduction and statement of results}

For $n \in \mathbb N$, $N > 0$, and $\alpha > -1/2$,
consider the unitary random matrix ensemble
\begin{equation}
\label{randommatrixmodel}
    Z_{n,N}^{-1} |\det M|^{2\alpha} e^{-N \Tr V(M)} \;dM,
\end{equation}
on the space $\mathcal{M}(n)$ of $n\times n$ Hermitian matrices $M$,
where $V$ is real analytic and satisfies
\begin{equation} \label{Vgrowth}
    \lim_{x \to \pm \infty} \frac{V(x)}{\log(x^2+1)} = +\infty.
    \end{equation}
Suppose also that the equilibrium measure $\mu_V$ for $V$ is regular
\cite{DKMVZ2}, and that $0$ is a right endpoint of the support of
$\mu_V$.

In the paper \cite{IKO1} we computed the limiting eigenvalue
correlation kernel in the double scaling limit as $n, N \to \infty$
such that $n^{2/3}(n/N-1) = O(1)$. We showed that it is
characterized through a solution of a model RH problem associated
with a special solution of the equation number XXXIV from the list
of Painlev\'e and Gambier \cite{Ince},
\begin{equation}
\label{painleve34}
    u'' = 4 u^2 + 2s u + \frac{(u')^2 - (2\alpha)^2}{2u}.
\end{equation}
In particular, we showed how the relevant solution $u(s)$ of
\eqref{painleve34} (which we denote by $u_{\alpha}(s)$)
can be obtained from a solution of a model RH
problem, which we next describe.

\subsection{The model RH problem} \label{modelRHP}
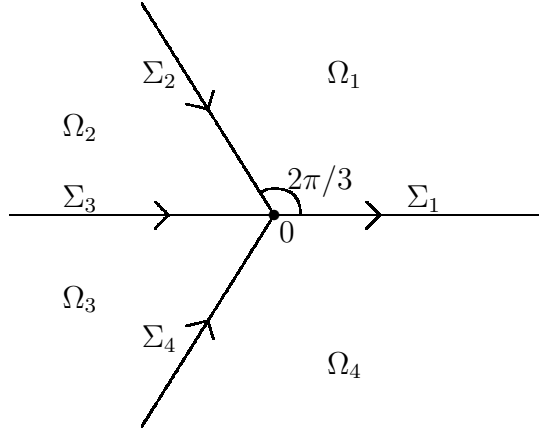
\begin{figure}[th]
\centering
\unitlength 1pt
\linethickness{0.5pt}
\begin{picture}(200,160)(-100,-80)
   \put(0,0){\line(1,0){100}}
   \put(0,0){\line(-1,0){100}}
   \qbezier(0,0)(-25,40)(-50,80)
   \qbezier(0,0)(-25,-40)(-50,-80)
   \qbezier(40,0)(39,1)(35,5)   \qbezier(40,0)(39,-1)(35,-5)
   \qbezier(-40,0)(-41,1)(-45,5)   \qbezier(-40,0)(-41,-1)(-45,-5)
   \qbezier(-25,40)(-24,43.5)(-23,47)   \qbezier(-25,40)(-29,41)(-33,42)
   \qbezier(-25,-40)(-24,-43.5)(-23,-47)   \qbezier(-25,-40)(-29,-41)(-33,-42)
   \put(2,-10){$0$}
   \put(0,0){\circle*{4}}
   \qbezier(10,0)(10,10)(0,10)
   \qbezier(0,10)(-4,10)(-5,8)
   \put(5,10){$2\pi/3$}
   \put(50,3){$\Sigma_1$}
   \put(-50,50){$\Sigma_{2}$}
   \put(-80,3){$\Sigma_3$}
   \put(-50,-50){$\Sigma_{4}$}
   \put(20,50){$\Omega_1$}
   \put(-80,30){$\Omega_2$}
   \put(-80,-35){$\Omega_3$}
   \put(20,-60){$\Omega_4$}
   \end{picture}
   \caption{Contours for the model Riemann-Hilbert problem.}\label{figure1}
\end{figure}

The model RH problem is posed on a contour $\Sigma$ in an auxiliary
$z$-plane, consisting of four rays $\Sigma_1 = \{ \arg z = 0\}$,
$\Sigma_2 = \{ \arg z = 2\pi/3\}$, $\Sigma_3 = \{ \arg z = \pi\}$,
and $\Sigma_4 = \{ \arg z = -2\pi/3\}$ with orientation as shown in
Figure~\ref{figure1}. As usual in RH problems, the orientation
defines a $+$ and a $-$ side on each part of the contour, where the
$+$-side is on the left when traversing the contour according to its
orientation. For a function $f$ on $\mathbb C \setminus \Sigma \equiv \Omega$, we
use $f_{\pm}$ to denote its limiting values on $\Sigma$ taken from
the $\pm$-side, provided such limiting values exist.
The contour $\Sigma$ divides the complex plane into four sectors $\Omega_j$ also
shown in the figure.

The model RH problem reads as follows.
\paragraph{Riemann-Hilbert problem for $\Psi_{\alpha}$}
\begin{enumerate}
\item[\rm (a)] $\Psi_{\alpha} : \mathbb{C} \setminus \Sigma  \to
    \mathbb C^{2\times 2}$ is analytic.
\item[\rm (b)] $\Psi_{\alpha,+}(z) = \Psi_{\alpha,-}(z)
    \begin{pmatrix} 1 & 1 \\ 0 & 1 \end{pmatrix}$,
    for $z \in \Sigma_1$,

    $\Psi_{\alpha,+}(z) = \Psi_{\alpha,-}(z)
    \begin{pmatrix} 1 & 0 \\ e^{2 \alpha \pi i} & 1 \end{pmatrix}$,
    for $z \in \Sigma_2$,

    $\Psi_{\alpha,+}(z) = \Psi_{\alpha,-}(z)
    \begin{pmatrix} 0 & 1 \\ -1 & 0  \end{pmatrix}$,
    for $z \in \Sigma_3$,

    $\Psi_{\alpha,+}(z) = \Psi_{\alpha,-}(z)
    \begin{pmatrix} 1 & 0 \\ e^{-2\alpha \pi i} & 1 \end{pmatrix}$,
    for $z \in \Sigma_4$.
\item[\rm (c)] $\Psi_{\alpha}(z) =  \left(I + O\left(\frac{1}{z}\right)\right)
z^{-\sigma_3/4}
    \frac{1}{\sqrt{2}} \begin{pmatrix} 1 & i \\ i & 1 \end{pmatrix}
e^{-(\frac{2}{3} z^{3/2} + s z^{1/2})\sigma_3}$
    as $z \to \infty$.

    Here $\sigma_3 = \left(\begin{smallmatrix} 1 & 0 \\ 0 & -1 \end{smallmatrix}\right)$
    is the third Pauli matrix.
\item[\rm (d)] $\Psi_{\alpha}(z) = O\begin{pmatrix} |z|^{\alpha} & |z|^{\alpha} \\
    |z|^{\alpha} & |z|^{\alpha}
\end{pmatrix}$ as $z \to 0$, if $-1/2 < \alpha < 0$; and

    $\Psi_{\alpha}(z) = \left\{ \begin{array}{ll}
    O\begin{pmatrix} |z|^{\alpha} & |z|^{-\alpha} \\
    |z|^{\alpha} & |z|^{-\alpha} \end{pmatrix}
    & \text{as $z \to 0$ with $z \in \Omega_1 \cup \Omega_4$}, \\[10pt]
    O\begin{pmatrix} |z|^{-\alpha} & |z|^{-\alpha} \\
    |z|^{-\alpha} & |z|^{-\alpha} \end{pmatrix}
    & \text{as $z \to 0$ with $z \in \Omega_2 \cup \Omega_3$},
    \end{array} \right.$
    if $\alpha \geq 0$.
 \end{enumerate}
   Here, and in what follows, the $O$-terms are taken entrywise. Note
that the RH problem depends on a parameter $s$ through the
asymptotic condition at infinity. If we want to emphasize the
dependence on $s$ we will write $\Psi_{\alpha}(z;s)$ instead of
$\Psi_{\alpha}(z)$.

The model RH problem is uniquely solvable for every $\alpha >
-1/2$ and $s \in \mathbb R$ (for details see \cite[Proposition 2.1]{IKO1}).

All solutions of \eqref{painleve34} are meromorphic in the complex
plane. The special solution of relevance in \cite{IKO1} is characterized
by the following result.
\begin{theorem} \label{theorem1}
Assume $\alpha > -1/2$.
Let $\Psi_{\alpha}(z;s)$ be the solution of the model RH problem
and write
\begin{equation} \label{Psiexpansion}
    \Psi_{\alpha}(z;s) =
    \left(I + \frac{m_{\Psi}(s)}{z} + O\left(\frac{1}{z^2}\right)\right)
z^{-\sigma_3/4}
    \frac{1}{\sqrt{2}} \begin{pmatrix} 1 & i \\ i & 1 \end{pmatrix}
e^{-(\frac{2}{3} z^{3/2} + s z^{1/2})\sigma_3}
\end{equation}
as $z \to \infty$.
Then
\begin{equation} \label{usolution}
    u_{\alpha}(s) = -\frac{s}{2} - i \frac{d}{ds}  \left[ m_{\Psi}(s) \right]_{12}
\end{equation}
exists and satisfies \eqref{painleve34}. The function
\eqref{usolution} is a global solution of \eqref{painleve34} (i.e.,
it does not have poles on the real line).

Moreover, $u_{\alpha}$ is also given by
\begin{equation} \label{usolutionat0}
    u_{\alpha}(s) = i \lim_{z \to 0} \left[ z
    \left(\frac{d}{dz} \Psi_{\alpha}(z)\right) \Psi_{\alpha}^{-1}(z) \right]_{12}.
    \end{equation}
\end{theorem}

\begin{proof}
 The expression \eqref{usolutionat0}
follows from Lemma 3.2 in \cite{IKO1}. The remaining statements of Theorem \ref{theorem1}
are contained in  Theorem 1.4 of \cite{IKO1}.
\end{proof}

\subsection{Main result} \label{main}
The aim of the present paper is an analysis of the asymptotic behavior of
the special solution $u_{\alpha}(s)$, described in Theorem \ref{theorem1},
as $s \rightarrow \pm \infty$.
Our main result is the following.
\begin{theorem} \label{theorem2}
Let $u_{\alpha}(s)$ be the solution of \eqref{painleve34} given in Theorem \ref{theorem1}.
Then,
\begin{equation}
\label{uplusinfinity}
u_{\alpha}(s) = \frac{\alpha}{\sqrt{s}} + O(s^{-2}),
\qquad \text{as } s \to +\infty,
\end{equation}
and
\begin{equation}
\label{uminusinfinity}
u_{\alpha}(s) = \frac{\alpha}{\sqrt{-s}}\cos \left(\frac{4}{3}(-s)^{3/2} -\alpha \pi\right)
+ O(s^{-2}), \qquad \text{as } s \to -\infty.
\end{equation}
\end{theorem}

In the Sections \ref{section2} -- \ref{section3}, we provide a proof
of our main result, Theorem \ref{theorem2}. This is accomplished by
using the Deift-Zhou steepest descent method for RH problems
\cite{DZ}. In the case at hand it consists of constructing a
sequence of invertible transformations $\Psi_{\alpha} \mapsto
A_{\alpha} \mapsto B_{\alpha} \mapsto C_{\alpha} \mapsto
D_{\alpha}$, where the matrix-valued function $D_{\alpha}$ is close
to the identity as $s \to \pm \infty$. By following the above
transformations asymptotics for $\Psi_{\alpha}$ and thus, in view of
\eqref{usolution} and \eqref{usolutionat0}, for $u_{\alpha}(s)$ may be derived.

\subsection{RH problem for Painlev\'e XXXIV} \label{RHproblem34}
In this paper we are mainly concerned with the special solution $u_{\alpha}(s)$.
The analysis of the general solution of the Painlev\'e XXXIV equation  \eqref{painleve34}
can be also performed via the nonlinear steepest descent method applied
to the following  generalization of the RH problem above.

\paragraph{Riemann-Hilbert problem for the general solution of PXXXIV}
\begin{enumerate}
\item[\rm (a)] $\Psi_{\alpha} : \mathbb{C} \setminus \Sigma  \to
    \mathbb C^{2\times 2}$ is analytic.
\item[\rm (b)] $\Psi_{\alpha,+}(z) = \Psi_{\alpha,-}(z)
    \begin{pmatrix} 1 & b_1 \\ 0 & 1 \end{pmatrix}$,
    for $z \in \Sigma_1$,

    $\Psi_{\alpha,+}(z) = \Psi_{\alpha,-}(z)
    \begin{pmatrix} 1 & 0 \\ b_2 & 1 \end{pmatrix}$,
    for $z \in \Sigma_2$,

    $\Psi_{\alpha,+}(z) = \Psi_{\alpha,-}(z)
    \begin{pmatrix} 0 & 1 \\ -1 & 0  \end{pmatrix}$,
    for $z \in \Sigma_3$,

    $\Psi_{\alpha,+}(z) = \Psi_{\alpha,-}(z)
    \begin{pmatrix} 1 & 0 \\ b_4 & 1 \end{pmatrix}$,
    for $z \in \Sigma_4$.
\item[\rm (c)] $\Psi_{\alpha}(z) =  \left(I + O\left(\frac{1}{z}\right)\right)
z^{-\sigma_3/4}
    \frac{1}{\sqrt{2}} \begin{pmatrix} 1 & i \\ i & 1 \end{pmatrix}
e^{-(\frac{2}{3} z^{3/2} + s z^{1/2})\sigma_3}$
    as $z \to \infty$.
\item[\rm (d)] If $\alpha - \frac{1}{2} \not\in \mathbb N_0$, then
\begin{equation}\label{connection0P34}
    \Psi_{\alpha}(z) = B(z)
    \begin{pmatrix} z^{\alpha} & 0 \\ 0 & z^{-\alpha} \end{pmatrix} E_j,
    \quad \text{ for } z \in \Omega_j,
\end{equation}
where $B$ is analytic. If $\alpha \in \frac{1}{2} + \mathbb N_0$, then there exists a constant $\kappa$ such that
\begin{equation} \label{connection1P34}
    \Psi_{\alpha}(z) = B(z)
    \begin{pmatrix} z^{\alpha} & \kappa z^{\alpha}\log z \\ 0 & z^{-\alpha} \end{pmatrix} E_j,
    \quad \text{ for } z \in \Omega_j,
\end{equation}
where $B$ is analytic.
\end{enumerate}
The complex numbers $b_{1}$, $b_{2}$, $b_{4}$ (the Stokes multipliers) and
the constant invertible matrices $E_{j}$, $j = 1, 2, 3, 4$ (the connection matrices)
form the RH data. They satisfy  certain general constraints which in particular
yield the following {\it cyclic} relation for the parameters $b_{j}$,
\begin{equation}\label{b_s}
b_{1} + b_{2} + b_{4} - b_{1}b_{2}b_{4} = 2\cos(2\alpha \pi).
\end{equation}
Except for the special case,
\begin{equation}\label{specP34}
 2\alpha =  n, \quad b_1 = b_2 = b_4 =  (-1)^{n},
 \quad n \in \mathbb N,
\end{equation}
when the solution of the RH problem is given in fact in terms of the Airy functions,
the connection matrices $E_{j}$ are determined
(up to inessential left diagonal or upper triangular factors) by $\alpha$ and the
Stokes multipliers $b_{j}$ \footnote{We refer to  \cite[Chapter 5]{FIKN} for more details on
the setting of the general RH problems for Painlev\'e equations.}.

In the formulation of the general RH problem we keep the previous notation $\Psi_{\alpha}(z)$
for its solution. Formulas \eqref{usolution} and \eqref{usolutionat0} for the solution of equation
\eqref{painleve34} are still valid, although, of course, the function $u_{\alpha}(s)$
for an arbitrary choice of the monodromy parameters $b_{j}$ might
have poles on the real line.

The case of our special interest in this paper corresponds to
the choice,
\begin{equation}\label{p34spec}
b_{1} = 1, \quad b_{2} = e^{2\alpha\pi i},
\quad b_{4} = e^{-2\alpha\pi i},
\end{equation}
of the Stokes parameters $b_{j}$. In Section \ref{section2} we
treat the case $s \to +\infty$, which
turns out to be the easier case.
In Section \ref{section3} we deal with the  more involved
case $s \to -\infty$. The main technical issue is the necessity
to construct an extra, in comparison with the
$+\infty$ case, parametrix  with Bessel functions.

It is well known (see e.g. \cite{Ince})  that the Painlev\'e XXXIV equation
\eqref{painleve34} can be in fact transformed to the Painlev\'e II equation.
We discuss in detail some aspects of this transformation, relevant to our analysis,
in the Appendix. In particular, we notice  that the second Painlev\'e
function which is associated with the special Painlev\'e XXXIV
solution $u_{\alpha}(s)$ we are studying here is {\it not}
the familiar in random matrix \cite{TW1}, \cite{TW2}, \cite{CKV} and string
\cite{SS} theories Hastings-McLeod function.
In addition, we show that, although the asymptotic behavior of
$u_{\alpha}(s)$ as $s \to +\infty$ can be extracted from the already known
asymptotics of the second Painlev\'e transcendents,
the asymptotic behavior  of  $u_{\alpha}(s)$ as $s \to -\infty$
needs indeed a separate analysis.

In the Appendix we also discuss the question of the uniqueness
of the solution $u_{\alpha}(s)$. In fact, we show that there is
a one-parameter family of solutions of equation
\eqref{painleve34} with the asymptotics \eqref{uplusinfinity}.
This family is characterizes by the choice,
\begin{equation}\label{family}
 b_{2} = e^{2\alpha\pi i},
\quad b_{4} = e^{-2\alpha\pi i},\quad
b_{1} = \text{arbitrary complex number},
\end{equation}
of the monodromy data $b_{j}$. At the same time, we
 conjecture that the
asymptotic condition \eqref{uminusinfinity} fixes the
solution uniquely. We present the arguments in support of this
conjecture which are based on certain observations related to the
asymptotic investigation of the general Painlev\'e XXXIV RH problem.

We conclude the introduction by mentioning that  the Painlev\'e XXXIV
equation has also appeared in several other physical applications.
In fact, in the already mentioned paper \cite{SS} it was the Painlev\'e XXXIV
image of the Hastings-McLeod Painlev\'e II solution which showed
up and not the solution itself.


\section{Proof of Theorem \ref{theorem2}: asymptotics as $s \to +\infty$}
\label{section2}

\subsection{First transformation $\Psi_{\alpha} \mapsto A_{\alpha}$}
Introduce
\begin{equation}
\label{re-scaling}
 A_{\alpha}(z;s) = s^{\sigma_3/4} \Psi_{\alpha}(sz;s), \qquad z \in \mathbb{C} \setminus \Sigma.
\end{equation}
It is then easy to see that $A_{\alpha}$ satisfies a RH problem similar to that
for $\Psi_{\alpha}$. In the following we often suppress the $s$-dependence
of functions whenever they are understood.
\paragraph{Riemann-Hilbert problem for $A_{\alpha}$}
\begin{enumerate}
\item[\rm (a)] $A_{\alpha} : \mathbb{C} \setminus \Sigma  \to
    \mathbb C^{2\times 2}$ is analytic.
\item[\rm (b)] $A_{\alpha}$ has the same jumps on $\Sigma$ as that of $\Psi_{\alpha}$.
\item[\rm (c)] $A_{\alpha}(z) = (I + O(1/z)) z^{-\sigma_3/4} \frac{1}{\sqrt{2}}
      \begin{pmatrix} 1 & i \\ i & 1 \end{pmatrix}
      e^{-s^{3/2}(\frac{2}{3} z^{3/2} + z^{1/2})\sigma_3}$
    as $z \to \infty$.
\item[\rm (d)] $A_{\alpha}$ has the same behavior near $0$ as that of $\Psi_{\alpha}$.
\end{enumerate}
From \eqref{usolutionat0} and \eqref{re-scaling} it follows that
\begin{equation} \label{usolution2}
    u_{\alpha}(s) = \frac{i}{\sqrt{s}} \lim_{z \to 0}  \left[ z
         \left(\frac{d}{dz} A_{\alpha}(z)\right) A_{\alpha}^{-1}(z) \right]_{12}.
\end{equation}

\subsection{Second transformation $A_{\alpha} \mapsto B_{\alpha}$}
\begin{figure}[th]
\centering
\unitlength 1pt
\linethickness{0.5pt}
\begin{picture}(200,160)(-100,-80)
   \put(0,0){\line(1,0){100}}
   \put(0,0){\line(-1,0){100}}
   \qbezier(0,0)(-25,40)(-50,80)
   \qbezier(0,0)(-25,-40)(-50,-80)
   \qbezier(-30,0)(-55,40)(-80,80)
   \qbezier(-30,0)(-55,-40)(-80,-80)
   \qbezier(40,0)(39,1)(35,5)   \qbezier(40,0)(39,-1)(35,-5)
   \qbezier(-10,0)(-11,1)(-15,5)   \qbezier(-10,0)(-11,-1)(-15,-5)
   \qbezier(-60,0)(-61,1)(-65,5)   \qbezier(-60,0)(-61,-1)(-65,-5)
   \qbezier(-25,40)(-24,43.5)(-23,47)   \qbezier(-25,40)(-29,41)(-33,42)
   \qbezier(-25,-40)(-24,-43.5)(-23,-47)   \qbezier(-25,-40)(-29,-41)(-33,-42)
   \qbezier(-55,40)(-54,43.5)(-53,47)   \qbezier(-55,40)(-59,41)(-63,42)
   \qbezier(-55,-40)(-54,-43.5)(-53,-47)   \qbezier(-55,-40)(-59,-41)(-63,-42)
   \put(2,-10){$0$}
   \put(0,0){\circle*{4}}
   \put(-35,-10){$-1$}
   \put(-30,0){\circle*{4}}
   \put(50,3){$\widetilde \Sigma_0$}
   \put(-28,3){$\widetilde \Sigma_1$}
   \put(-25,50){$\Sigma_{2}$}
   \put(-75,45){$\widetilde \Sigma_2$}
   \put(-90,3){$\widetilde \Sigma_3$}
   \put(-25,-55){$\Sigma_{4}$}
    \put(-75,-50){$\widetilde \Sigma_4$}
   \put(20,50){I}
   \put(-50,50){II}
   \put(-90,30){III}
   \put(-90,-30){IV}
   \put(-50,-50){V}
   \put(20,-60){VI}
   \end{picture}
   \caption{Contours and domains for the definition of $B_{\alpha}$.}\label{figure2}
\end{figure}
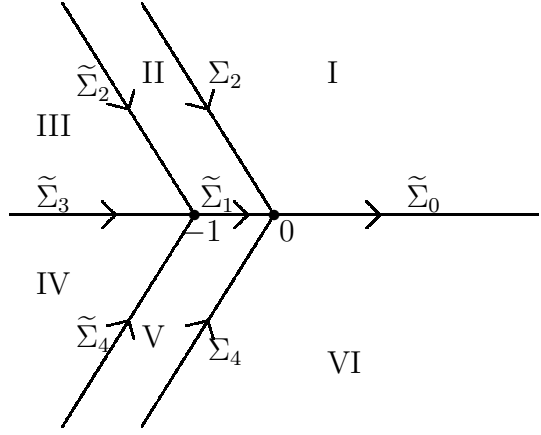

Let us put
\begin{equation}\label{tposit}
t = s^{3/2},
\end{equation}
which is the large parameter in the RH problem.
Consider Figure~\ref{figure2}.

Define
\begin{equation} \label{Bdef}
    B_{\alpha}(z) = \left\{ \begin{array}{ll}
        A_{\alpha}(z), & \quad \text{ for } z \in I \cup III \cup IV \cup VI,\\
        A_{\alpha}(z)\begin{pmatrix} 1 & 0 \\ e^{2\alpha \pi i} & 1\end{pmatrix},
        & \quad \text{ for } z \in II,\\[10pt]
        A_{\alpha}(z)\begin{pmatrix} 1 & 0 \\ - e^{-2\alpha \pi i} & 1\end{pmatrix},
        & \quad \text{ for } z \in V.
        \end{array} \right.
    \end{equation}
Let $\widetilde{\Sigma} = \Sigma - 1$, i.e., $\widetilde{\Sigma}$ is
$\Sigma$ translated to the left by one. Then
\begin{equation}
\widetilde{\Sigma} = \bigcup_{j=0}^4 \widetilde{\Sigma}_j,
\end{equation}
where the disjoint contours $\widetilde{\Sigma}_j$ are oriented as in Figure~\ref{figure2}.

\paragraph{Riemann-Hilbert problem for $B_{\alpha}$}
\begin{enumerate}
\item[\rm (a)] $B_{\alpha} : \mathbb{C} \setminus \widetilde{\Sigma}  \to
    \mathbb C^{2\times 2}$ is analytic.
\item[\rm (b)] $B_{\alpha,+}(z) = B_{\alpha,-}(z)
    \begin{pmatrix} 1 & 1 \\ 0 & 1 \end{pmatrix}$,
    for $z \in \widetilde{\Sigma}_0 = (0,\infty)$,

    $B_{\alpha,+}(z) = B_{\alpha,-}(z)
    \begin{pmatrix} e^{2\alpha \pi i} & 1 \\ 0 & e^{-2\alpha \pi i}\end{pmatrix}$,
    for $z \in \widetilde{\Sigma}_1 = (-1,0)$,

    $B_{\alpha,+}(z) = B_{\alpha,-}(z)
    \begin{pmatrix} 1 & 0 \\ e^{2 \alpha \pi i} & 1 \end{pmatrix}$,
    for $z \in \widetilde{\Sigma}_2$,

    $B_{\alpha,+}(z) = B_{\alpha,-}(z)
    \begin{pmatrix} 0 & 1 \\ -1 & 0  \end{pmatrix}$,
    for $z \in \widetilde{\Sigma}_3 = (-\infty,-1)$,

    $B_{\alpha,+}(z) = B_{\alpha,-}(z)
    \begin{pmatrix} 1 & 0 \\ e^{-2\alpha \pi i} & 1 \end{pmatrix}$,
    for $z \in \widetilde{\Sigma}_4$.
\item[\rm (c)] $B_{\alpha}(z) = (I + O(1/z)) z^{-\sigma_3/4}
    \frac{1}{\sqrt{2}} \begin{pmatrix} 1 & i \\ i & 1 \end{pmatrix}
     e^{-t(\frac{2}{3} z^{3/2} + z^{1/2})\sigma_3}$
    as $z \to \infty$.
\item[\rm (d)] If $\alpha < 0$, then $B_{\alpha}(z) = O\begin{pmatrix} |z|^{\alpha} & |z|^{\alpha} \\
    |z|^{\alpha} & |z|^{\alpha} \end{pmatrix}$ as $z \to 0$.\\[5pt]
    If $\alpha \geq 0$, then $B_{\alpha}(z) =
    O\begin{pmatrix} |z|^{\alpha} & |z|^{-\alpha} \\
    |z|^{\alpha} & |z|^{-\alpha} \end{pmatrix}$ as $z \to 0$.
\end{enumerate}

From \eqref{usolution2} and \eqref{Bdef} it follows that
\begin{equation} \label{usolution3}
    u_{\alpha}(s) = \frac{i}{\sqrt{s}} \lim_{z \to 0}  \left[ z
         \left(\frac{d}{dz} B_{\alpha}(z)\right) B_{\alpha}^{-1}(z) \right]_{12}.
\end{equation}

\subsection{Third transformation $B_{\alpha} \mapsto C_{\alpha}$}
We next introduce the $g$-function
\begin{equation}\label{gposit}
g(z) = \frac{2}{3} (z+1)^{3/2},\qquad -\pi < \arg (z+1) < \pi.
\end{equation}
By a straightforward computation
\begin{equation}
\label{gasympt}
g(z) - \Big(\frac{2}{3}z^{3/2} + z^{1/2}\Big) = \frac{1}{4} z^{-1/2} + O(z^{-3/2}),
\qquad \text{as } z \to \infty.
\end{equation}
Define
\begin{equation} \label{Calphadef}
C_{\alpha}(z) = \begin{pmatrix} 1 & 0 \\ - i t/4 & 1 \end{pmatrix}
    B_{\alpha}(z) e^{t g(z) \sigma_3}.
\end{equation}
Then $C_{\alpha}$ satisfies the following RH problem.
\paragraph{Riemann-Hilbert problem for $C_{\alpha}$}
\begin{enumerate}
\item[\rm (a)] $C_{\alpha} : \mathbb{C} \setminus \widetilde{\Sigma} \to
    \mathbb C^{2\times 2}$ is analytic.
\item[\rm (b)] $C_{\alpha,+}(z) = C_{\alpha,-}(z)
    \begin{pmatrix} 1 & e^{-2t g(z)} \\ 0 & 1 \end{pmatrix}$,
    for $z \in \widetilde{\Sigma}_0$,

    $C_{\alpha,+}(z) = C_{\alpha,-}(z)
    \begin{pmatrix} e^{2\alpha \pi i} & e^{-2t g(z)} \\ 0 & e^{-2\alpha \pi i}\end{pmatrix}$,
    for $z \in \widetilde{\Sigma}_1$,

    $C_{\alpha,+}(z) = C_{\alpha,-}(z)
    \begin{pmatrix} 1 & 0 \\ e^{2 \alpha \pi i} e^{2t g(z)} & 1 \end{pmatrix}$,
    for $z \in \widetilde{\Sigma}_2$,

    $C_{\alpha,+}(z) = C_{\alpha,-}(z)
    \begin{pmatrix} 0 & 1 \\ -1 & 0  \end{pmatrix}$,
    for $z \in \widetilde{\Sigma}_3$,

    $C_{\alpha,+}(z) = C_{\alpha,-}(z)
    \begin{pmatrix} 1 & 0 \\ e^{-2\alpha \pi i} e^{2t g(z)} & 1 \end{pmatrix}$,
    for $z \in \widetilde{\Sigma}_4$.
\item[\rm (c)] $C_{\alpha}(z) = (I + O(1/z)) z^{-\sigma_3/4}
    \frac{1}{\sqrt{2}} \begin{pmatrix} 1 & i \\ i & 1 \end{pmatrix}$
    as $z \to \infty$.
\item[\rm (d)] If $\alpha < 0$, then $C_{\alpha}(z) = O\begin{pmatrix} |z|^{\alpha} & |z|^{\alpha} \\
    |z|^{\alpha} & |z|^{\alpha} \end{pmatrix}$ as $z \to 0$.\\[5pt]
    If $\alpha \geq 0$, then $C_{\alpha}(z) =
    O\begin{pmatrix} |z|^{\alpha} & |z|^{-\alpha} \\
    |z|^{\alpha} & |z|^{-\alpha} \end{pmatrix}$ as $z \to 0$.
\end{enumerate}
Note that the prefactor $\left(\begin{smallmatrix} 1 & 0 \\ - i t/4 & 1 \end{smallmatrix}\right)$
in the definition (\ref{Calphadef}) of $C_{\alpha}$ is needed for the asymptotic condition (c) in the RH
problem. The prefactor does not affect the $12$ entry and so does
not influence the computation of $u_{\alpha}$ via the formula \eqref{usolution3}.

Thus by \eqref{Calphadef}
\begin{equation} \label{Calphacalculation}
\left[ z \left(\frac{d}{dz} B_{\alpha}(z) \right) B_{\alpha}^{-1}(z) \right]_{12}
    =
    \left[ z \left(\frac{d}{dz} C_{\alpha}(z) \right) C_{\alpha}^{-1}( z) \right]_{12}
        - tg'(z) \left[ z C_{\alpha}(z) \sigma_3 C_{\alpha}^{-1}(z) \right]_{12}.
        \end{equation}
In view of item (d) in the RH problem for $C_{\alpha}$ (and the fact that
$\det C_{\alpha} \equiv 1$) we have that
$z C_{\alpha}(z) \sigma_3 C_{\alpha}^{-1}(z) \to 0$ as $ z\to 0$.
Therefore the second term in the right-hand side of \eqref{Calphacalculation}
vanishes as $z \to 0$ and it follows by \eqref{usolution3} that
\begin{equation} \label{usolution4}
    u_{\alpha}(s) = \frac{i}{\sqrt{s}} \lim_{z \to 0}  \left[ z
         \left(\frac{d}{dz} C_{\alpha}(z)\right) C_{\alpha}^{-1}(z) \right]_{12}.
\end{equation}

\subsection{Construction of parametrices}
\subsubsection{Global parametrix $P_{\alpha}^{(\infty)}$}
Away from the point $-1$ we expect that $C_{\alpha}$ should be well approximated
by the solution $P^{(\infty)}_{\alpha}$ of the following RH problem.
\paragraph{Riemann-Hilbert problem for $P^{(\infty)}_{\alpha}$}
\begin{enumerate}
\item[\rm (a)] $P^{(\infty)}_{\alpha} : \mathbb{C} \setminus (-\infty, 0] \to
    \mathbb C^{2\times 2}$ is analytic.
\item[\rm (b)] $P^{(\infty)}_{\alpha,+}(z) = P^{(\infty)}_{\alpha,-}(z)
    \begin{pmatrix} 0 & 1 \\ -1 & 0\end{pmatrix}$,
    for $z \in (-\infty,-1)$,

    $P^{(\infty)}_{\alpha,+}(z) = P^{(\infty)}_{\alpha,-}(z)
    \begin{pmatrix} e^{2\alpha \pi i} & 0\\ 0 & e^{-2\alpha \pi i} \end{pmatrix}$,
    for $z \in (-1,0)$.

\item[\rm (c)] $P^{(\infty)}_{\alpha}(z) =
    (I + O(1/z)) z^{-\sigma_3/4}
    \frac{1}{\sqrt{2}} \begin{pmatrix} 1 & i \\ i & 1 \end{pmatrix}$
    as $z \to \infty$.
\end{enumerate}
It should be noted that $P^{(\infty)}_{\alpha}$ does not depend on $s$.

We seek $P^{(\infty)}_{\alpha}$ in the form
\begin{equation}
P^{(\infty)}_{\alpha}(z) = F_{\alpha}(z) z^{\alpha \sigma_3},
\end{equation}
where $F_{\alpha}$ is analytic in $\mathbb{C} \setminus (-\infty,-1]$. Clearly
then the jump is correct on $(-1,0)$. A straightforward computation
shows that in order to have the correct jump also on $(-\infty, -1)$ we may take
\begin{equation}
F_{\alpha}(z) = E (z+1)^{-\sigma_3/4} \frac{1}{\sqrt{2}}
\begin{pmatrix} 1 & i \\ i & 1\end{pmatrix} (\delta_{\alpha}(z))^{\sigma_3},
\end{equation}
where $E$ is a constant prefactor and
\begin{equation}
\delta_{\alpha}(z) = \exp \left(-\frac{\alpha}{\pi} (z+1)^{1/2}
\int_{1}^{\infty} \frac{\log t}{\sqrt{t-1}(t+z)}\,dt \right).
\end{equation}
Using the residue theorem and a contour deformation argument, it is also straightforward
to see that
\begin{alignat}{2}
\nonumber \int_{1}^{\infty} \frac{\log t}{\sqrt{t-1}(t+z)}\,dt &=
\frac{\pi \log z}{(z+1)^{1/2}} + \pi \int_0^1 \frac{dt}{\sqrt{1-t} (t+z)}\\[5pt]
&= \frac{\pi \log z}{(z+1)^{1/2}} +
\frac{\pi\,\log \Big( \frac{(z + 1)^{1/2} + 1}{(z + 1)^{1/2} - 1} \Big)}{(z+1)^{1/2}}.
\end{alignat}
Hence,
\begin{alignat}{2}
\nonumber
P^{(\infty)}_{\alpha}(z) & = E
    (z+1)^{-\sigma_3/4} \frac{1}{\sqrt{2}}
    \begin{pmatrix} 1 & i \\ i & 1\end{pmatrix}
     (\delta_{\alpha}(z))^{\sigma_3} z^{\alpha \sigma_3}\\
& = E (z+1)^{-\sigma_3/4} \frac{1}{\sqrt{2}} \begin{pmatrix} 1 & i \\ i & 1\end{pmatrix}
    \bigg( \frac{(z + 1)^{1/2} + 1}{(z + 1)^{1/2} - 1} \bigg)^{-\alpha \sigma_3}.
\label{pinfp}
\end{alignat}
In order to satisfy the asymptotic condition (c) of the RH problem we should take
\begin{align} \label{defE}
    E = \begin{pmatrix} 1 & 0 \\ 2\alpha i & 1 \end{pmatrix}.
\end{align}

Note that
\[ \frac{d}{dz} \log \left(\frac{(z + 1)^{1/2} + 1}{(z + 1)^{1/2} - 1}\right) =
    - \frac{1}{z (z+1)^{1/2}} \]
from which it follows after straightforward calculations from \eqref{pinfp}
and \eqref{defE} that
\begin{align} \nonumber
    \lim_{ z \to 0} z
    \left(\frac{d}{dz} P^{(\infty)}_{\alpha}(z) \right)
    \left(P^{(\infty)}_{\alpha}(z) \right)^{-1}
    & = E \frac{1}{\sqrt{2}} \begin{pmatrix} 1 & i \\ i & 1 \end{pmatrix}
        \left(\alpha \sigma_3 \right)
        \frac{1}{\sqrt{2}} \begin{pmatrix} 1 & -i \\ -i & 1 \end{pmatrix} E^{-1} \\
        & \label{Pinftylogderivative}
        = - i \alpha
            \begin{pmatrix} - 2 \alpha i & 1 \\
                -1 + 4\alpha^2 & 2 \alpha i \end{pmatrix}.
\end{align}

\subsubsection{Local parametrix $P^{(-1)}_{\alpha}$}
The global parametrix $P^{(\infty)}_{\alpha}$ will not be a good approximation
to $C_{\alpha}$ near the point $-1$.
Let $U^{(-1)}$ be a small open disc around $-1$ of radius $< 1$.
We seek a local parametrix $P^{(-1)}_{\alpha}$
defined in $U^{(-1)}$ which satisfies the following.

\paragraph{Riemann-Hilbert problem for $P^{(-1)}_{\alpha}$}
\begin{enumerate}
\item[(a)] $P^{(-1)}_{\alpha} : \overline{U^{(-1)}} \setminus \widetilde{\Sigma} \to
    \mathbb C^{2\times 2}$ is continuous and analytic on $U^{(-1)} \setminus
    \widetilde{\Sigma}$.
\item[(b)] $P^{(-1)}_{\alpha,+}(z) = P^{(-1)}_{\alpha,-}(z)\,v_{C_{\alpha}}(z)$ for $z \in
    \widetilde{\Sigma} \cap U^{(-1)}$, where $v_{C_{\alpha}}$ denotes the jump matrix
    for $C_{\alpha}$
   (the contour having the same orientation as $\widetilde{\Sigma}$).
\item[(c)] $P^{(-1)}_{\alpha}(z)\,\left(P^{(\infty)}_{\alpha}(z)\right)^{-1} = I +
O\big(\frac{1}{t}\big)$, as $t \to \infty$,
uniformly for $z \in \partial U^{(-1)} \setminus \widetilde{\Sigma}$.
\end{enumerate}

We seek $P^{(-1)}_{\alpha}$ in the form
\begin{equation}\label{P1n}
P^{(-1)}_{\alpha}(z) = \widehat{P}^{(-1)}_{\alpha}(z) e^{t g(z) \sigma_3},
\end{equation}
where $\widehat{P}^{(-1)}_{\alpha}$ satisfies the following RH problem with constant jumps.

\paragraph{Riemann-Hilbert problem for $\widehat{P}^{(-1)}_{\alpha}$}
\begin{enumerate}
\item[(a)] $\widehat{P}^{(-1)}_{\alpha} : \overline{U^{(-1)}} \setminus \widetilde{\Sigma} \to
    \mathbb C^{2\times 2}$ is continuous and analytic on $U^{(-1)} \setminus \Sigma_S$.
\item[(b)]
$\widehat{P}^{(-1)}_{\alpha,+}(z) = \widehat{P}^{(-1)}_{\alpha,-}(z)
    \begin{pmatrix}
        e^{2\alpha \pi i} & 1 \\
        0 & e^{-2\alpha \pi i}
    \end{pmatrix}$, for $z \in \widetilde{\Sigma}_1 \cap U^{(-1)}$,

$\widehat{P}^{(-1)}_{\alpha,+}(z) = \widehat{P}^{(-1)}_{\alpha,-}(z)
    \begin{pmatrix}
        1 & 0 \\
        e^{2 \alpha \pi i} & 1
    \end{pmatrix}$, for $z \in \widetilde{\Sigma}_2 \cap U^{(-1)}$,

$\widehat{P}^{(-1)}_{\alpha,+}(z) = \widehat{P}^{(-1)}_{\alpha,-}(z)
    \begin{pmatrix}
        0 & 1 \\
        -1 & 0
    \end{pmatrix}$, for $z \in \widetilde{\Sigma}_3 \cap U^{(-1)}$,

$\widehat{P}^{(-1)}_{\alpha,+}(z) = \widehat{P}^{(-1)}_{\alpha,-}(z)
    \begin{pmatrix}
        1 & 0 \\
        e^{-2 \alpha \pi i} & 1
    \end{pmatrix}$, for $z \in \widetilde{\Sigma}_4 \cap U^{(-1)}$.
\item[(c)] $\widehat{P}^{(-1)}_{\alpha}(z) =
  P^{(\infty)}_{\alpha}(z) \Big( I + O\big(\frac{1}{t}\big)\Big) e^{-tg(z) \sigma_3}$, as $t \to \infty$,
uniformly for $z \in \partial U^{(-1)} \setminus \widetilde{\Sigma}$.
\end{enumerate}

A solution to this RH problem can be constructed in terms of Airy functions.
The standard Airy parametrix is posed in an auxiliary $\zeta$-plane
and satisfies the following RH problem for a contour $\Sigma$ as in Figure~\ref{figure1}.
\paragraph{Riemann-Hilbert problem for $\Phi^{(Ai)}$}
\begin{enumerate}
\item[(a)] $\Phi^{(Ai)} : \mathbb C \setminus \Sigma \to \mathbb C^{2 \times 2}$ is analytic.
\item[(b)] $\Phi^{(Ai)}_+ = \Phi^{(Ai)}_-
    \begin{pmatrix}
        1 & 1 \\
        0 & 1
    \end{pmatrix}$, on $\Sigma_1$,

$\Phi^{(Ai)}_+ = \Phi^{(Ai)}_-
    \begin{pmatrix}
        1 & 0 \\
        1 & 1
    \end{pmatrix}$, on $\Sigma_2 \cup \Sigma_4$,

    $\Phi^{(Ai)}_+ = \Phi^{(Ai)}_-
    \begin{pmatrix}
        0 & 1 \\
        -1 & 0
    \end{pmatrix}$, on $\Sigma_3$.
\item[(c)] $\Phi^{(Ai)}(\zeta) = \zeta^{-\sigma_3/4} (I + O(\zeta^{-3/2}))
    \frac{1}{2\sqrt{\pi}} \begin{pmatrix} 1 & i \\ -1 & i \end{pmatrix}
     e^{-\frac{2}{3} \zeta^{3/2} \sigma_3}$
        as $\zeta \to \infty$.
\end{enumerate}
The solution is built out of the functions
\[ y_0(\zeta) = \Ai(\zeta), \quad y_1(\zeta) = \omega \Ai(\omega \zeta), \quad
 y_2(\zeta) = \omega^2 \Ai(\omega^2 \zeta), \qquad \omega = e^{2\pi i/3}, \]
and takes the following form
\begin{equation}
    \left\{ \begin{array}{llcll}
    \Phi^{(Ai)}  = \begin{pmatrix} -y_1 & -y_2 \\ -y_1' & -y_2' \end{pmatrix}
    & \text{in } \Omega_2, & \quad &
    \Phi^{(Ai)} = \begin{pmatrix} y_0 & -y_2 \\ y_0' & -y_2' \end{pmatrix}
    & \text{in } \Omega_1, \\[10pt]
    \Phi^{(Ai)}  = \begin{pmatrix} -y_2 & y_1 \\ -y_2' & y_1' \end{pmatrix}
    & \text{in } \Omega_3, & \quad &
    \Phi^{(Ai)}  = \begin{pmatrix} y_0 & y_1 \\ y_0' & y_1' \end{pmatrix}
    & \text{in } \Omega_4.
    \end{array} \right.
\end{equation}
Then we put
\begin{equation}  \label{Phiadef}
    \Phi^{(Ai)}_{\alpha}(\zeta) = \sqrt{2\pi} \begin{pmatrix} 1 & 0 \\
    0 & -i \end{pmatrix} \Phi^{(Ai)}(\zeta) e^{\pm \alpha \pi i \sigma_3}, \qquad
        \text{for } \pm \Im \zeta > 0,
        \end{equation}
and $\Phi^{(Ai)}_{\alpha}$ satisfies the following RH problem.
\paragraph{Riemann-Hilbert problem for $\Phi^{(Ai)}_{\alpha}$}
\begin{enumerate}
\item[(a)] $\Phi^{(Ai)}_{\alpha} : \mathbb{C} \setminus \Sigma \to
    \mathbb C^{2\times 2}$ is analytic.
\item[(b)]
$\Phi^{(Ai)}_{\alpha,+} = \Phi^{(Ai)}_{\alpha,-}
    \begin{pmatrix}
        e^{2\alpha \pi i} & 1 \\
        0 & e^{-2\alpha \pi i}
    \end{pmatrix}$, on $\Sigma_1$,

$\Phi^{(Ai)}_{\alpha,+} = \Phi^{(Ai)}_{\alpha,-}
    \begin{pmatrix}
        1 & 0 \\
        e^{2 \alpha \pi i} & 1
    \end{pmatrix}$,  on $\Sigma_2$,

$\Phi^{(Ai)}_{\alpha,+} = \Phi^{(Ai)}_{\alpha,-}
    \begin{pmatrix}
        0 & 1 \\
        -1 & 0
    \end{pmatrix}$, on $\Sigma_3$,

$\Phi^{(Ai)}_{\alpha,+} = \Phi^{(Ai)}_{\alpha,-}    \begin{pmatrix}
        1 & 0 \\
        e^{-2 \alpha \pi i} & 1
    \end{pmatrix}$, on $\Sigma_4$.
\item[(c)] $\ds \Phi^{(Ai)}_{\alpha}(\zeta) =
    \zeta^{-\sigma_3/4}\Big( I + O(\zeta^{-3/2}) \Big)
        \frac{1}{\sqrt{2}} \begin{pmatrix} 1 & i \\ i & 1 \end{pmatrix}
            e^{\pm \alpha \pi i \sigma_3}
    e^{-\frac{2}{3} \zeta^{3/2} \sigma_3}$

        as $\zeta \to \infty$ with $\pm \Im \zeta > 0$.
\end{enumerate}

Define
\begin{equation}\label{P1n2}
\widehat{P}^{(-1)}_{\alpha}(z) =  E_{\alpha}(z) \Phi^{(Ai)}_{\alpha}(s(z+1)),
\quad \text{for } z \in U^{(-1)} \setminus \widetilde{\Sigma},
\end{equation}
where $E_{\alpha}$ is analytic in $U^{(-1)}$. Then $\widehat{P}^{(-1)}_{\alpha}$ has
the correct jumps. In order to satisfy the matching condition
in the RH problem we take $E_{\alpha}$ in the following way:
\begin{equation}\label{Eadef}
E_{\alpha}(z) = P_{\alpha}^{(\infty)}(z)  e^{\mp \alpha \pi i \sigma_3}
    \frac{1}{\sqrt{2}} \begin{pmatrix} 1 & -i \\ -i & 1 \end{pmatrix}
    (s(z+1))^{\sigma_3/4}.
\end{equation}
It is a straightforward computation to verify that $E_{\alpha}$ extends
as an analytic function in $U^{(-1)}$. Combining (\ref{P1n}),
(\ref{P1n2}), and (\ref{Eadef}), we see that
\[ P_{\alpha}^{(-1)}(z) = P_{\alpha}^{(\infty)}(z)  e^{\mp \alpha \pi i \sigma_3}
    \frac{1}{\sqrt{2}} \begin{pmatrix} 1 & -i \\ -i & 1 \end{pmatrix}
    (s(z+1))^{\sigma_3/4} \Phi^{(Ai)}_{\alpha}(s(z+1)) e^{tg(z)\sigma_3}, \]
which completes the construction
of the local parametrix $P_{\alpha}^{(-1)}$.

\subsection{Fourth transformation $C_{\alpha} \mapsto D_{\alpha}$}
Define now the final transformation
\begin{equation} \label{Ddef}
    D_{\alpha}(z) = \left\{ \begin{array}{ll}
        C_{\alpha}(z) \big(P^{(-1)}_{\alpha}(z)\big)^{-1}, & \quad \text{ for }
        z \in U^{(-1)} \setminus \widetilde{\Sigma},\\[10pt]
        C_{\alpha}(z) \big(P^{(\infty)}_{\alpha}(z)\big)^{-1}, & \quad \text{ for }
        z \in \mathbb{C} \setminus (\overline{U^{(-1)} \cup \widetilde{\Sigma}}).
        \end{array} \right.
\end{equation}

\begin{figure}[th]
\centering
\unitlength 1pt
\linethickness{0.5pt}
\begin{picture}(200,160)(-100,-80)
   \put(0,0){\line(1,0){100}}
   \put(0,0){\line(-1,0){10}}
   \qbezier(-40.5,16.8)(-55,40)(-80,80)
   \qbezier(-40.5,-16.8)(-55,-40)(-80,-80)
   \put(-30,0){\circle{40}}
   \qbezier(40,0)(39,1)(35,5)   \qbezier(40,0)(39,-1)(35,-5)
   \qbezier(-50,0)(-51,1)(-55,5)   \qbezier(-50,0)(-49,1)(-45,5)
   \qbezier(-55,40)(-54,43.5)(-53,47)   \qbezier(-55,40)(-59,41)(-63,42)
   \qbezier(-55,-40)(-54,-43.5)(-53,-47)   \qbezier(-55,-40)(-59,-41)(-63,-42)
   \put(2,-10){$0$}
   \put(0,0){\circle*{4}}
   \put(-35,-10){$-1$}
   \put(-30,0){\circle*{4}}
   \end{picture}
   \caption{Contour $\Sigma_{D}$ in the RH problem for $D_{\alpha}$.}\label{figure3}
\end{figure}
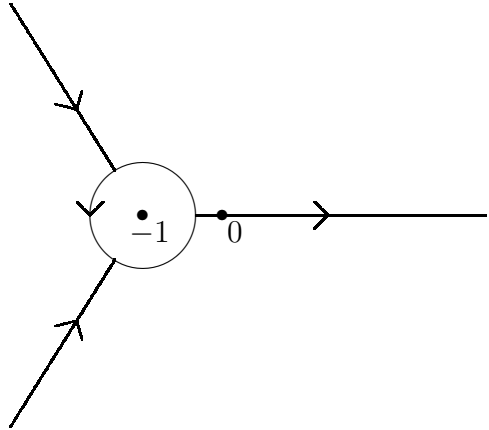

Since $C_{\alpha}$ and $P_{\alpha}^{(\infty)}$ have the same jumps on
$(-\infty, -1) \setminus U^{(-1)}$, and $C_{\alpha}$ and $P_{\alpha}^{(-1)}$
have the same jumps on
$U^{(-1)} \cap \widetilde{\Sigma}$,
we have that $D_{\alpha}$ is analytic across these contours. What remains are
jumps for $D_{\alpha}$ on the contour $\Sigma_{D}$ shown in Figure~\ref{figure3}.
Indeed, $D_{\alpha}$ satisfies the following RH problem.

\paragraph{Riemann-Hilbert problem for $D_{\alpha}$}
\begin{enumerate}
\item[(a)] $D_{\alpha} : \mathbb{C} \setminus \Sigma_{D} \to \mathbb C^{2\times 2}$ is analytic.
\item[(b)] $D_{\alpha,+}(z) = D_{\alpha,-}(z)\,v_{D_{\alpha}}(z)$ for $z \in \Sigma_{D}$, where\\
\begin{equation} \label{vRdef}
v_{D_{\alpha}} = \left\{
    \begin{array}{ll}
        P^{(\infty)}_{\alpha}\,(P^{(-1)}_{\alpha})^{-1},&\text{ on } \partial U^{(-1)},\\
        P^{(\infty)}_{\alpha, -}\,v_{C_{\alpha}}\,(P^{(\infty)}_{\alpha,+})^{-1},&\text{ on }
        \Sigma_{D} \setminus \partial U^{(-1)}.
\end{array}
\right.
\end{equation}
\item[(c)] $D_{\alpha}(z) = I + O(1/z)$ as $z \to \infty$.
\end{enumerate}
The subscripts $\pm$ in  $P^{(\infty)}$ are only relevant for the
segment of the horizontal part of the contour to the left of $0$.

\subsection{Conclusion of the proof of \eqref{uplusinfinity}}
The jump matrix satisfies
\begin{equation}\label{vD1}
v_{D_{\alpha}}(z) = I + O(1/t), \quad \text{ as } t \to \infty,
\end{equation}
uniformly on the circle $\partial U^{(-1)}$. In addition
\begin{equation}\label{vD2}
v_{D_{\alpha}}(z) = I + O(e^{-ct(|z|+1)}), \quad c > 0,
\end{equation}
uniformly on $\Sigma_{D} \setminus \partial U^{(-1)}$.
In a standard way (see e.g. \cite{IKO1}) one shows that
\begin{equation}
    \label{Dalphaasympt1}
D_{\alpha}(z) = I + O\left(\frac{1}{t(1+|z|)}\right), \quad \text{ as } t \to \infty,
\end{equation}
uniformly for $z \in \mathbb{C} \setminus \Sigma_{D}$.

Finally, we have by  \eqref{Ddef} and the fact that $D_{\alpha}(z)$
and $\frac{d}{dz} D_{\alpha}(z)$ remain bounded as $z \to 0$,
\[ \lim_{z \to \infty} z \left(\frac{d}{dz} C_{\alpha}(z) \right)  C_{\alpha}^{-1}(z)
    = D_{\alpha}(0) \lim_{z \to \infty} z \left( \frac{d}{dz} P_{\alpha}^{(\infty)}(z) \right)
    \left( P_{\alpha}^{(\infty)}(z) \right)^{-1} D_{\alpha}^{-1}(0)
    \]
so that in view of
\eqref{usolution4} and
\eqref{Pinftylogderivative}
\begin{equation} \label{usolution5}
    u_{\alpha}(s) = \frac{\alpha}{\sqrt{s}}
    \left[ D_{\alpha}(0)
        \begin{pmatrix} - 2\alpha i & 1 \\ -1 + 4\alpha^2 & 2 \alpha i \end{pmatrix}
            D_{\alpha}^{-1}(0) \right]_{12}.
            \end{equation}
Inserting \eqref{Dalphaasympt1} with $z = 0$ into \eqref{usolution5}
and recalling that $t = s^{3/2}$, we
obtain \eqref{uplusinfinity}.

\section{Proof of Theorem \ref{theorem2}: asymptotics as $s \to -\infty$}
\label{section3}

For the asymptotics as $s \to -\infty$ we also perform a sequence
of transformations of the model RH problem
$\Psi_{\alpha} \mapsto A_{\alpha} \mapsto B_{\alpha} \mapsto C_{\alpha} \mapsto D_{\alpha}$,
but the transformations are different from the ones we performed
for $s \to +\infty$. Thus $A_{\alpha}$,
$B_{\alpha}$, $C_{\alpha}$ and $D_{\alpha}$ will now have a different meaning
which hopefully does not lead to any confusion.
 We assume throughout this section that $s < 0$.

\subsection{First transformation $\Psi_{\alpha} \mapsto A_{\alpha}$}
Similar to (\ref{re-scaling}), we introduce
\begin{equation}
\label{re-scaling2}
A_{\alpha}(z;s) = (-s)^{\sigma_{3}/4}
\Psi_{\alpha}(-sz;s), \qquad z \in \mathbb{C} \setminus \Sigma.
\end{equation}
The $A_{\alpha}$ - RH problem  reads as follows
(the $s$-dependence is, as usual, suppressed).
\paragraph{Riemann-Hilbert problem for $A_{\alpha}$}
\begin{enumerate}
\item[\rm (a)] $A_{\alpha} : \mathbb{C} \setminus \Sigma  \to
    \mathbb C^{2\times 2}$ is analytic.
\item[\rm (b)] The jumps of $A_{\alpha}$ on $\Sigma$ are the same
as those of $\Psi_{\alpha}$.
\item[\rm (c)] $A_{\alpha}(z) = \left(I + O\left(\frac{1}{z}\right)\right)
z^{-\sigma_3/4}
    \frac{1}{\sqrt{2}} \begin{pmatrix} 1 & i \\ i & 1 \end{pmatrix}
e^{-t(\frac{2}{3} z^{3/2} - z^{1/2})\sigma_3}$
as $z \to \infty$.
\item[\rm (d)] $A_{\alpha}$ has the same behavior near $0$ as $\Psi_{\alpha}$ has.
\end{enumerate}
Here, the large positive parameter $t$ is defined by the equation
(cf. \eqref{tposit})
 \begin{equation} \label{tneg}
t = (-s)^{3/2}.
\end{equation}

Using \eqref{usolutionat0} together with \eqref{re-scaling2} we can
express $u_{\alpha}$ in terms of $A_{\alpha}$ as follows
\begin{equation} \label{usolution2A}
    u_{\alpha}(s) = \frac{i}{\sqrt{-s}} \lim_{z \to 0}  \left[ z
         \left(\frac{d}{dz} A_{\alpha}(z)\right) A_{\alpha}^{-1}(z) \right]_{12}.
\end{equation}

\subsection{Second transformation $A_{\alpha} \mapsto B_{\alpha}$}
An important difference comparing with the previous case is
that a step analogous to the $B_{\alpha}$ - step is skipped. That is, our
next step will be the $g$-function ``dressing''.

Put (cf.(\ref{gposit}))
\begin{equation} \label{gneg}
g(z) = \frac{2}{3} (z-1)^{3/2},\quad -\pi < \arg (z-1) < \pi.
\end{equation}
Note that, as before,
\begin{equation}
\label{gasymptneg}
g(z) - \Big(\frac{2}{3}z^{3/2} - z^{1/2}\Big) = \frac{1}{4} z^{-1/2} + O(z^{-3/2}),
\qquad \text{as } z \to \infty.
\end{equation}
Define
\begin{equation} \label{fromAtoB}
B_{\alpha}(z) =
\begin{pmatrix} 1 & 0\\ -\frac{it}{4} & 1 \end{pmatrix}
A_{\alpha}(z) e^{t g(z) \sigma_3}.
\end{equation}
Then, $B_{\alpha}$ satisfies the following RH problem.
\paragraph{Riemann-Hilbert problem for $B_{\alpha}$}
\begin{enumerate}
\item[\rm (a)] $B_{\alpha} : \mathbb{C} \setminus \Sigma \to
    \mathbb C^{2\times 2}$ is analytic.
\item[\rm (b)] $B_{\alpha,+}(z) = B_{\alpha,-}(z)
    \begin{pmatrix} 1 & e^{-2t g(z)} \\ 0 & 1 \end{pmatrix}$,
    for $z \in (1,\infty)$,

    $B_{\alpha,+}(z) = B_{\alpha,-}(z)
    \begin{pmatrix} e^{-t(g_{-}(z) - g_{+}(z))} & 1
    \\ 0 &  e^{t(g_{-}(z) - g_{+}(z))}\end{pmatrix}$,
    for $z \in (0,1)$,

    $B_{\alpha,+}(z) = B_{\alpha,-}(z)
    \begin{pmatrix} 1 & 0 \\ e^{2 \alpha \pi i + 2t g(z)} & 1 \end{pmatrix}$,
    for $z \in \Sigma_2$,

    $B_{\alpha,+}(z) = B_{\alpha,-}(z)
    \begin{pmatrix} 0 & 1 \\ -1 & 0  \end{pmatrix}$,
    for $z \in \Sigma_3$,

    $B_{\alpha,+}(z) = B_{\alpha,-}(z)
    \begin{pmatrix} 1 & 0 \\ e^{-2\alpha \pi i + 2t g(z)} & 1 \end{pmatrix}$,
    for $z \in \Sigma_4$.
\item[\rm (c)] $B_{\alpha}(z) =
 \left(I + O\left(\frac{1}{z}\right)\right)
z^{-\sigma_3/4}
    \frac{1}{\sqrt{2}} \begin{pmatrix} 1 & i \\ i & 1 \end{pmatrix}$
as $z \to \infty$.
\item[\rm (d)] The behavior of $B_{\alpha}(z)$ as $z \to 0$
is the same as that of $\Psi_{\alpha}(z)$.
\end{enumerate}
We emphasize, that the contour $\Sigma$ is now {\it the same} as
in the original $\Psi_{\alpha}$ - problem.

From the transformation \eqref{fromAtoB} it follows that
\begin{equation} \label{logderivativeB}
    \left[\frac{d}{dz} A_{\alpha}(z) A_{\alpha}^{-1}(z) \right]_{12}
    = \left[\frac{d}{dz} B_{\alpha}(z) B_{\alpha}^{-1}(z) \right]_{12}
    - t g'(z) \left[B_{\alpha}(z) \sigma_3 B_{\alpha}^{-1}(z) \right]_{12}.
    \end{equation}
From part (d) in the RH problem satisfied by $B_{\alpha}$ we can
deduce that
\[ B_{\alpha}(z) \sigma_3 B_{\alpha}^{-1}(z) =
    \begin{cases} O(|z|^{2\alpha}) & \text{ if } -1/2 < \alpha < 0, \\
    O(1) & \text{ if } \alpha \geq 0,
    \end{cases} \]
as $z \to 0$.
It follows that we can forget about the second term in the right-hand
side of \eqref{logderivativeB} and we obtain from \eqref{usolution2A}
and \eqref{logderivativeB} that
\begin{equation} \label{usolution3B}
    u_{\alpha}(s) = \frac{i}{\sqrt{-s}}
        \lim_{z \to 0} \left[ z \left(\frac{d}{dz} B_{\alpha}(z) \right) B_{\alpha}^{-1}(z) \right]_{12}.
\end{equation}

\subsection{Third transformation $B_{\alpha} \mapsto C_{\alpha}$}

\begin{figure}[htb]
\centering \begin{overpic}[scale=0.6]%
{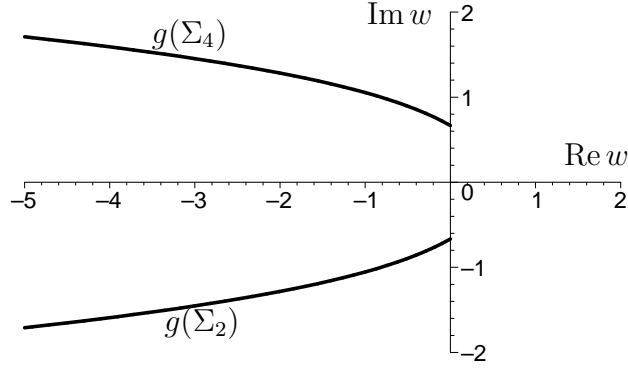}%
       \put(25,5){$g(\Sigma_2)$}
       \put(23,52){$g(\Sigma_4)$}
       \put(90,32){$\Re w$}
       \put(58,55){$\Im w$}
\end{overpic}
\caption{Images of the rays $\Sigma_2$ and $\Sigma_4$ under the mapping
$z \mapsto w = g(z) = \frac{2}{3}(z-1)^{3/2}$.}\label{figure4}
\end{figure}

In order to proceed further, we need to analyze the structure of the
sign of $\Re g(z)$. The first observation is trivial:
\begin{equation}\label{g1}
g(z) = \frac{2}{3}|z-1|^{3/2} \geq 0, \quad z \in (1,\infty).
\end{equation}
Next, we notice that the function $ w = g(z)$ performs a conformal
mapping of the upper half plane to the sector $0 < \arg w < 3\pi/2$.
Under this mapping, the domain $ \frac{\pi}{3} < \arg (z-1) < \pi$
becomes the left half plane  $ \Re w <0$, and the ray $\Sigma_{2}$
transforms to a simple smooth curve $g(\Sigma_2)$ which
ends at $-2i/3$, lies entirely in the left half
plane, and behaves for large $z \in \Sigma_2$ as
\[ \Re w \sim - \frac{2}{3} |z|^{3/2},
    \qquad \Im w \sim -\frac{1}{2} \sqrt{3} |z|^{1/2}, \]
see Figure~\ref{figure4}.
Similarly, the function $ w = g(z)$  performs a conformal
mapping of the lower half plane to the sector $-3\pi/2 < \arg w < 0$.
Under this mapping, the domain $-\pi < \arg (z-1) < - \frac{\pi}{3} $
becomes the left half plane  $ \Re w <0$, and the ray $\Sigma_{4}$
is mapped to the mirror image of  $g(\Sigma_2)$ with respect to the
real axis, see again Figure~\ref{figure4}.
Therefore, there exists a constant $c > 0$ such that
\begin{equation}\label{g2}
\Re g(z)\leq - c|z-1|\leq  0, \quad z \in \Sigma_{2}\cup \Sigma_{4}.
\end{equation}

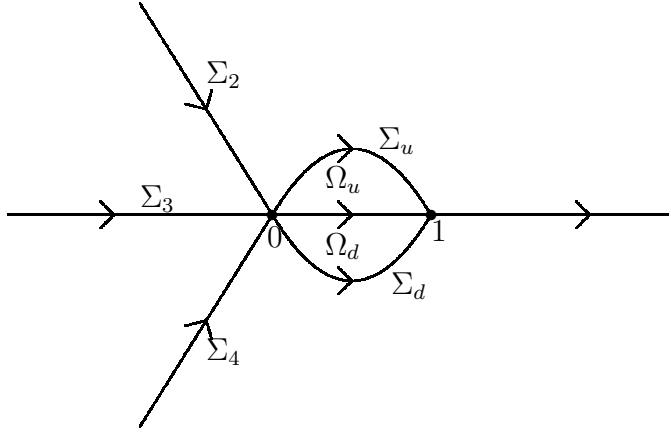
\begin{figure}[th]
\centering
\unitlength 1pt
\linethickness{0.5pt}
\begin{picture}(200,160)(-60,-80)
   \put(0,0){\line(1,0){150}}
   \put(0,0){\line(-1,0){100}}
   \qbezier(0,0)(-25,40)(-50,80)
   \qbezier(0,0)(-25,-40)(-50,-80)
   \qbezier(0,0)(30,50)(60,0)
   \qbezier(0,0)(30,-50)(60,0)
   \qbezier(120,0)(119,1)(115,5)   \qbezier(120,0)(119,-1)(115,-5)
   \qbezier(30,0)(29,1)(25,5)   \qbezier(30,0)(29,-1)(25,-5)
   \qbezier(-60,0)(-61,1)(-65,5)   \qbezier(-60,0)(-61,-1)(-65,-5)
   \qbezier(-25,40)(-24,43.5)(-23,47)   \qbezier(-25,40)(-29,41)(-33,42)
   \qbezier(-25,-40)(-24,-43.5)(-23,-47)   \qbezier(-25,-40)(-29,-41)(-33,-42)
   \qbezier(30,25)(29,26)(25,30) \qbezier(30,25)(29,24)(25,20)
   \qbezier(30,-25)(29,-26)(25,-30) \qbezier(30,-25)(29,-24)(25,-20)
   \put(-2,-12){$0$}
   \put(0,0){\circle*{4}}
   \put(60,-10){$1$}
   \put(60,0){\circle*{4}}
   \put(-25,50){$\Sigma_{2}$}
   \put(-50,3){$\Sigma_3$}
   \put(-25,-55){$\Sigma_{4}$}
   \put(40,25){$\Sigma_u$}
   \put(45,-30){$\Sigma_d$}
    \put(20,10){$\Omega_u$}
    \put(20,-15){$\Omega_d$}
   \end{picture}
   \caption{Contour $\Sigma_C = \Sigma \cup \Sigma_u \cup \Sigma_d$
   and domains  for the definition of $C_{\alpha}$.}\label{figure5}
\end{figure}

We also notice that the function,
\begin{equation}
h(z):= g_{-}(z) - g_{+}(z),\quad z \in (0,1),
\end{equation}
admits  analytic continuation into the domains
$\Omega_{u}$ and $\Omega_{d}$ indicated in Figure~\ref{figure5}.
Indeed we have,
\begin{equation}\label{hup}
h(z) = -2g(z), \quad z\in \Omega_{u},
\end{equation}
and
\begin{equation}\label{hdown}
h(z) = 2g(z), \quad z\in \Omega_{d}.
\end{equation}
The indicated above characterization of the conformal
mapping generated by the function $g(z)$ yields the
inequalities
\begin{equation}\label{g3}
\Re h(z) > 0, \quad z\in \Omega_{u},
\end{equation}
and
\begin{equation}\label{g4}
\Re h(z) <0, \quad z\in \Omega_{d}.
\end{equation}

The estimates (\ref{g1}) and (\ref{g2}) imply that the
jump matrices on the rays $(1,\infty)$,
$\Sigma_{2}$ and $\Sigma_{4}$ and away of the end points
are close to the identity matrix, while the
estimates (\ref{g3}) and (\ref{g4}) suggest to
``open the lenses''  around the interval  $(0,1)$.

Noticing that
\begin{align} \nonumber
  \begin{pmatrix} e^{-t(g_{-}(z) - g_{+}(z))} & 1
    \\ 0 &  e^{t(g_{-}(z) - g_{+}(z))}\end{pmatrix}
    & =
      \begin{pmatrix} e^{-th(z)} & 1
    \\ 0 &  e^{th(z)}\end{pmatrix} \\
    & = \label{factorization}
  \begin{pmatrix} 1 & 0
    \\ e^{th(z)} &  1\end{pmatrix}
   \begin{pmatrix} 0& 1
    \\ -1 &  0\end{pmatrix}
 \begin{pmatrix} 1 & 0
    \\ e^{-th(z)} &  1\end{pmatrix},
\end{align}
we define the new function $C_{\alpha}$ with the
help of the following equations.
\begin{equation} \label{Ydef}
    C_{\alpha}(z) = \left\{ \begin{array}{ll}
        B_{\alpha}(z), & \quad \text{ for } z \notin \Omega_{u}\cup \Omega_{d} ,\\
        B_{\alpha}(z)\begin{pmatrix} 1 & 0 \\- e^{-th(z)} & 1\end{pmatrix},
        & \quad \text{ for } z \in \Omega_{u},\\[10pt]
        B_{\alpha}(z)\begin{pmatrix} 1 & 0 \\  e^{th(z)} & 1\end{pmatrix},
        & \quad \text{ for } z \in \Omega_{d}.
        \end{array} \right.
    \end{equation}
We use $\Sigma_{u}$ and   $\Sigma_{d}$ to denote the curves which,
in conjunction  with the interval $[0,1]$,  make
the boundary of the lenses $\Omega_{u}$ and $\Omega_{d}$,
respectively. The curves  $\Sigma_{u}$ and   $\Sigma_{d}$ together
with their orientation are indicated in Figure~\ref{figure5}.
Let $\Sigma_C$ denote the contour $\Sigma$ augmented by the arcs
$\Sigma_{u}$ and   $\Sigma_{d}$.
Then, $C_{\alpha}$ satisfies the following RH problem.
\paragraph{Riemann-Hilbert problem for $C_{\alpha}$}
\begin{enumerate}
\item[\rm (a)] $C_{\alpha} : \mathbb{C} \setminus \Sigma_C \to
    \mathbb C^{2\times 2}$ is analytic.
\item[\rm (b)] $C_{\alpha,+}(z) = C_{\alpha,-}(z)
    \begin{pmatrix} 1 & e^{-2t g(z)} \\ 0 & 1 \end{pmatrix}$,
    for $z \in (1,\infty)$,

    $C_{\alpha,+}(z) = C_{\alpha,-}(z)
    \begin{pmatrix} 0& 1
    \\ -1 & 0\end{pmatrix}$,
    for $z \in (0,1) \cup \Sigma_3$,

    $C_{\alpha,+}(z) = C_{\alpha,-}(z)
    \begin{pmatrix} 1 & 0 \\  e^{2tg(z)} & 1 \end{pmatrix}$,
    for $z \in \Sigma_u \cup \Sigma_d$,

    $C_{\alpha,+}(z) = C_{\alpha,-}(z)
    \begin{pmatrix} 1 & 0 \\ e^{2\alpha \pi i + 2t g(z)} & 1 \end{pmatrix}$,
    for $z \in \Sigma_2$,

    $C_{\alpha,+}(z) = C_{\alpha,-}(z)
    \begin{pmatrix} 1 & 0 \\ e^{-2\alpha \pi i+2t g(z)} & 1 \end{pmatrix}$,
    for $z \in \Sigma_4$.
\item[\rm (c)] $C_{\alpha}(z) =
\left(I + O\left(\frac{1}{z}\right)\right)
z^{-\sigma_3/4}
    \frac{1}{\sqrt{2}} \begin{pmatrix} 1 & i \\ i & 1 \end{pmatrix}$
    as $z \to \infty$.
\item[\rm (d)]
$C_{\alpha}(z) = O\begin{pmatrix} |z|^{\alpha} & |z|^{\alpha} \\
    |z|^{\alpha} & |z|^{\alpha}
\end{pmatrix}$ as $z \to 0$, if $-1/2 < \alpha < 0$; and

    $C_{\alpha}(z) = \left\{ \begin{array}{ll}
    O\begin{pmatrix} |z|^{\alpha} & |z|^{-\alpha} \\
    |z|^{\alpha} & |z|^{-\alpha} \end{pmatrix}
    & \text{as $z \to 0$ with $z \in (\Omega_1 \cup \Omega_4) \setminus (\Omega_u \cup \Omega_d)$}, \\[10pt]
    O\begin{pmatrix} |z|^{-\alpha} & |z|^{-\alpha} \\
    |z|^{-\alpha} & |z|^{-\alpha} \end{pmatrix}
    & \text{as $z \to 0$ with $z \in \Omega_2 \cup \Omega_3 \cup \Omega_u \cup \Omega_d$},
    \end{array} \right.$
    if $\alpha \geq 0$.
\end{enumerate}
When formulating the jump conditions across the lenses boundaries,
i.e., on the curves $\Sigma_{u}$ and   $\Sigma_{d}$, we have replaced
the function $h(z)$ by the function $g(z)$ according to the
relations (\ref{hup}) and  (\ref{hdown}).
The contours for the $C_{\alpha}$ - RH problem are depicted in Figure~\ref{figure5}.

To express $u_{\alpha}$ in terms of $C_{\alpha}$ we can use
the same formula \eqref{usolution3B} but with the understanding
that $z \to 0$ from outside the lens. Thus
\begin{equation} \label{usolution4C}
    u_{\alpha}(s) = \frac{i}{\sqrt{-s}}
        \lim_{z \to 0 \atop z \not\in \Omega_u \cup \Omega_d}
        \left[ z \left(\frac{d}{dz} C_{\alpha}(z) \right) C_{\alpha}^{-1}(z) \right]_{12}.
\end{equation}

\subsection{Construction of parametrices}

\subsubsection{Global parametrix $P^{(\infty)}$}
Away from the points $0$ and $1$ the jump matrices on $(1,\infty)$, $\Sigma_2$, $\Sigma_4$,
$\Sigma_u$, and $\Sigma_d$ all tend to the identity matrix as $t \to +\infty$
at an exponential rate. Therefore, away
from the points $0$ and $1$ we expect that $C_{\alpha}$ should be well
approximated by the solution $P^{(\infty)}$ of the
following RH problem with the only nontrivial jump across
the ray $(-\infty, 1)$.

\paragraph{Riemann-Hilbert problem for $P^{(\infty)}$}
\begin{enumerate}
\item[\rm (a)] $P^{(\infty)} : \mathbb{C} \setminus (-\infty, 1] \to
    \mathbb C^{2\times 2}$ is analytic.
\item[\rm (b)] $P^{(\infty)}_{+}(z) = P^{(\infty)}_{-}(z)
    \begin{pmatrix} 0 & 1 \\ -1 & 0\end{pmatrix}$,
    for $z \in (-\infty,1)$.
\item[\rm (c)] $P^{(\infty)}(z) =
\left(I + O\left(\frac{1}{z}\right)\right)
z^{-\sigma_3/4}
    \frac{1}{\sqrt{2}} \begin{pmatrix} 1 & i \\ i & 1 \end{pmatrix}$
    as $z \to \infty$.
\end{enumerate}
It should be noted that now  $P^{(\infty)}$ is clearly independent
of both $s$ and $\alpha$.

This Riemann-Hilbert problem is even simpler than the corresponding
problem for the case of positive $s$, and its solution is obviously
given by the formula (cf. (\ref{pinfp}))
\begin{equation}\label{Pinfsol}
P^{(\infty)}(z) = (z-1)^{-\sigma_3/4} \frac{1}{\sqrt{2}}
    \begin{pmatrix} 1 & i \\ i & 1\end{pmatrix}, \quad -\pi < \arg (z-1) < \pi.
\end{equation}

Near the points $0$ and $1$ the parametrix $P^{(\infty)}$ cannot be expected to
represent the asymptotics of  $C_{\alpha}$.
Our next task
is to construct the parametrix solutions near the points mentioned.

\subsubsection{Local parametrix $P^{(1)}$}
We begin
with the point $z = 1$. We will see that the form of parametrix
at  this point is very similar to the parametrix
near the point $z = -1$ in the previous case of positive $s$.

Let $U^{(1)}$ be a small open disc around $1$. We seek a parametrix $P^{(1)}$
defined in $U^{(1)}$ which satisfies the following RH problem.

\paragraph{Riemann-Hilbert problem for $P^{(1)}$}
\begin{enumerate}
\item[(a)] $P^{(1)} : \overline{U^{(1)}} \setminus \Sigma_C \to
    \mathbb C^{2\times 2}$ is continuous and analytic on $U^{(1)} \setminus
   \Sigma_C$.
\item[(b)] $P^{(1)}_{+}(z) = P^{(1)}_{-}(z)\,v_{C_{\alpha}}(z)$ for $z \in
    \Sigma_C \cap U^{(1)}$.
\item[(c)] $P^{(1)}(z)\,\left(P^{(\infty)}(z)\right)^{-1} = I +
O\big(\frac{1}{t}\big)$ as $t \to \infty$,
uniformly for $z \in \partial U^{(1)} \setminus \Sigma_C$.
\end{enumerate}
Similar to the case of positive $s$, we use  $v_{C_{\alpha}}$ to denote
the jump matrix in the RH problem for $C_{\alpha}$.
Comparing this Riemann-Hilbert problem with the one
for the function $P^{(-1)}_{\alpha}$ from the previous section,
we see that the solution $P^{(1)}$ can be given in terms of
the  matrix function $\Phi^{(Ai)}_{\alpha}$ defined by
equation (\ref{Phiadef}) and evaluated at $\alpha =0$. Indeed,
we propose the following form for the function
 $P^{(1)}$ (cf. (\ref{P1n}), (\ref{P1n2}))  :
\begin{equation}\label{P1def}
P^{(1)}(z) =  (-s)^{\sigma_{3}/4}\Phi^{(Ai)}_{0}( -s(z-1))e^{t g(z) \sigma_3},
\quad \text{for } z \in U^{(1)} \setminus \Sigma_C,
\end{equation}
where $\Phi^{(Ai)}_{0} :=\Phi^{(Ai)}_{\alpha = 0}$. We note that the relevant left multiplier
 $E(z)$ is chosen to be a scaling factor $(-s)^{\sigma_3/4}$, and that this is
exactly what the matrix function $E_{\alpha}(z)$ from \eqref{Eadef}
reduces to if $\alpha =0$ and if we replace
$s(z+1)$ by $-s(z-1)$ as is appropriate in the present situation.

Defined by (\ref{P1def}), the function $P^{(1)}(z)$ has obviously the
correct jumps. Moreover, when $z$ belongs to the
boundary of $U^{(1)}$ and $-s$ is large, the function
$\Phi_{0}(-s(z-1))$ in the r.h.s. of (\ref{P1def}) can be replaced by
its asymptotics. Therefore, we have:
\begin{align} \nonumber
P^{(1)}(z) & = (z-1)^{-\sigma_{3}/4}
\frac{1}{\sqrt{2}} \begin{pmatrix} 1 & i \\ i & 1\end{pmatrix}
\left(I + O\left(\frac{1}{(-s)^{3/2}}\right)\right) \\
& = \label{P1as1}
P^{(\infty)}(z)\Bigl(I + O\left(\frac{1}{t}\right)\Bigr),
\end{align}
as $ t \to \infty$,  uniformly for $
z \in \partial U^{(1)} \setminus \Sigma_C$,
which yields the matching condition needed. The
parametrix in the neighborhood of the point $z=1$
is then constructed. Note that it does not depend on $\alpha$.

\subsection{Local parametrix $P^{(0)}_{\alpha}$}

We are passing now to the analysis of the
neighborhood of the point $z=0$. Let $U^{(0)} $ be a small open disc around $0$.
We seek a parametrix $P^{(0)}_{\alpha}$ defined in $U^{(0)}$
which satisfies the following RH problem.

\paragraph{Riemann-Hilbert problem for $P^{(0)}_{\alpha}$}
\begin{enumerate}
\item[(a)] $P^{(0)}_{\alpha} : \overline{U^{(0)}} \setminus \Sigma_C \to
    \mathbb C^{2\times 2}$ is continuous and analytic on $U^{(0)} \setminus
   \Sigma_C$.
\item[(b)] $P^{(0)}_{\alpha,+}(z) = P^{(0)}_{\alpha,-}(z)\,v_{C_{\alpha}}(z)$ for $z \in
    \Sigma_C \cap U^{(0)}$.
\item[(c)] $P^{(0)}_{\alpha}(z)\left(P^{(\infty)}(z)\right)^{-1} = I +
O\big(\frac{1}{t}\big)$ as $t \to \infty$,
uniformly for $z \in \partial U^{(0)} \setminus \Sigma_C$.
\item[(d)] $P^{(0)}_{\alpha}$ has the same behavior near $0$ as $C_{\alpha}$.
\end{enumerate}

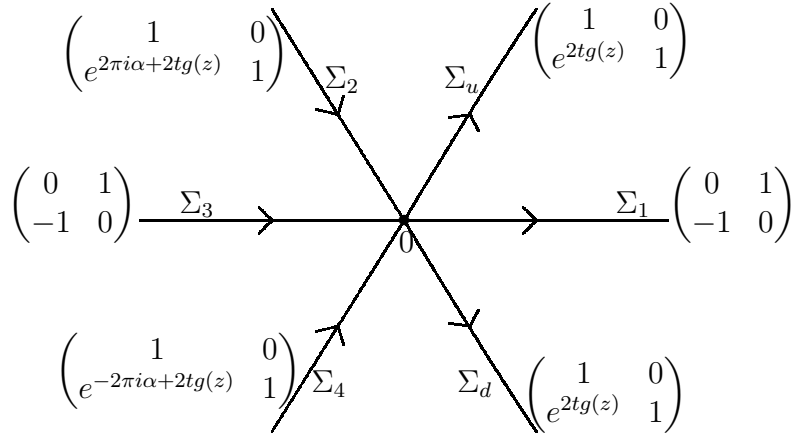
\begin{figure}[th]
\centering
\unitlength 1pt
\linethickness{0.5pt}
\begin{picture}(200,160)(-100,-80)
   \put(0,0){\line(1,0){100}}
   \put(0,0){\line(-1,0){100}}
   \qbezier(0,0)(-25,40)(-50,80)
   \qbezier(0,0)(-25,-40)(-50,-80)
   \qbezier(0,0)(30,50)(50,80)
   \qbezier(0,0)(30,-50)(50,-80)
   \qbezier(50,0)(49,1)(45,5)   \qbezier(50,0)(49,-1)(45,-5)
   \qbezier(-50,0)(-51,1)(-55,5)   \qbezier(-50,0)(-51,-1)(-55,-5)
   \qbezier(-25,40)(-24,43.5)(-23,47)   \qbezier(-25,40)(-29,41)(-33,42)
   \qbezier(-25,-40)(-24,-43.5)(-23,-47)   \qbezier(-25,-40)(-29,-41)(-33,-42)
   \qbezier(25,40)(21,39.5)(17,39)   \qbezier(25,40)(26.5,37)(28,34)
   \qbezier(25,-40)(21,-39.5)(17,-39)   \qbezier(25,-40)(26.5,-37)(28,-34)
   \put(-2,-12){$0$}
   \put(0,0){\circle*{4}}
   \put(80,3){$\Sigma_1$}
    \put(-30,50){$\Sigma_2$}
    \put(-85,3){$\Sigma_3$}
    \put(-35,-65){$\Sigma_4$}
    \put(15,50){$\Sigma_u$}
    \put(20,-65){$\Sigma_d$}
   \put(100,3){$\begin{pmatrix} 0& 1  \\ -1 & 0\end{pmatrix}$}
   \put(-130,60){$\begin{pmatrix} 1 & 0 \\  e^{2\pi i\alpha + 2tg(z)} & 1 \end{pmatrix}$}
   \put(-150,3){$\begin{pmatrix} 0& 1  \\ -1 & 0\end{pmatrix}$}
   \put(-132,-60){$\begin{pmatrix} 1 & 0 \\ e^{-2\pi i\alpha + 2tg(z)} & 1 \end{pmatrix}$}
   \put(47,65){$\begin{pmatrix} 1 & 0 \\  e^{2tg(z)} & 1 \end{pmatrix}$}
   \put(45,-69){$\begin{pmatrix} 1 & 0 \\  e^{2tg(z)} & 1 \end{pmatrix}$}
   \end{picture}
   \caption{Contours and jump matrices for the RH problem for $P_{\alpha}^{(0)}$
   (magnified picture).}\label{figure6}
\end{figure}

The contours and the jump matrices in the RH problem for $P_{\alpha}^{(0)}$
are depicted in Figure~\ref{figure6}.

\begin{figure}[th]
\centering
\unitlength 1pt
\linethickness{0.5pt}
\begin{picture}(200,160)(-100,-80)
   \put(0,0){\line(1,0){100}}
   \put(0,0){\line(-1,0){100}}
   \qbezier(0,0)(-25,40)(-50,80)
   \qbezier(0,0)(-25,-40)(-50,-80)
   \qbezier(0,0)(30,50)(50,80)
   \qbezier(0,0)(30,-50)(50,-80)
   \qbezier(50,0)(49,1)(45,5)   \qbezier(50,0)(49,-1)(45,-5)
   \qbezier(-50,0)(-51,1)(-55,5)   \qbezier(-50,0)(-51,-1)(-55,-5)
   \qbezier(-25,40)(-24,43.5)(-23,47)   \qbezier(-25,40)(-29,41)(-33,42)
   \qbezier(-25,-40)(-24,-43.5)(-23,-47)   \qbezier(-25,-40)(-29,-41)(-33,-42)
   \qbezier(25,40)(21,39.5)(17,39)   \qbezier(25,40)(26.5,37)(28,34)
   \qbezier(25,-40)(21,-39.5)(17,-39)   \qbezier(25,-40)(26.5,-37)(28,-34)
   \put(-2,-12){$0$}
   \put(0,0){\circle*{4}}
   \put(80,3){$\Sigma_1$}
    \put(-30,50){$\Sigma_2$}
    \put(-85,3){$\Sigma_3$}
    \put(-35,-65){$\Sigma_4$}
    \put(15,50){$\Sigma_u$}
    \put(20,-65){$\Sigma_d$}
   \put(100,3){$\begin{pmatrix} 0& 1  \\ -1 & 0\end{pmatrix}$}
   \put(-100,60){$\begin{pmatrix} 1 & 0 \\  e^{2\pi i\alpha} & 1 \end{pmatrix}$}
   \put(-150,3){$\begin{pmatrix} 0& 1  \\ -1 & 0\end{pmatrix}$}
   \put(-102,-60){$\begin{pmatrix} 1 & 0 \\ e^{-2\pi i\alpha} & 1 \end{pmatrix}$}
   \put(47,65){$\begin{pmatrix} 1 & 0 \\  1 & 1 \end{pmatrix}$}
   \put(45,-69){$\begin{pmatrix} 1 & 0 \\ 1 & 1 \end{pmatrix}$}

   \end{picture}
   \caption{Contours and jump matrices for the RH problem for $\widehat{P}_{\alpha}^{(0)}$
   (magnified picture).}\label{figure7}
\end{figure}
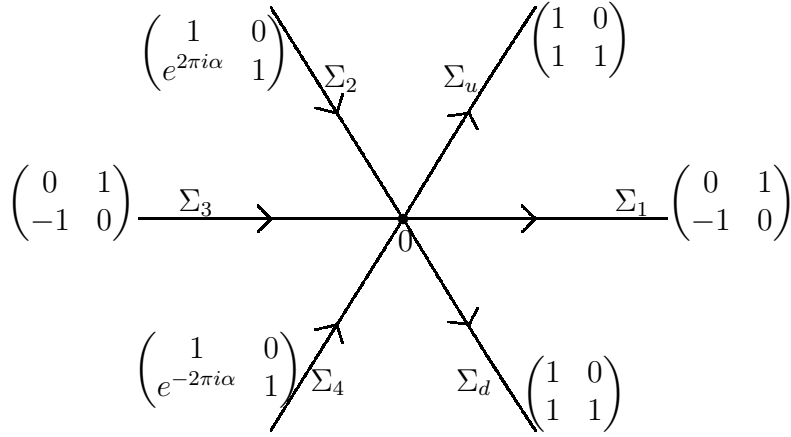

We take $P_{\alpha}^{(0)}$ in the form
\begin{equation} \label{defPalpha0}
    P_{\alpha}^{(0)}(z) = \widehat{P}_{\alpha}^{(0)}(z)
    e^{tg(z)\sigma_3}
    \end{equation}
and we see that $\widehat{P}_{\alpha}^{(0)}$ should
satisfy  a RH problem with jumps that are indicated in Figure~\ref{figure7}.
The jump matrices for $\widehat{P}_{\alpha}^{(0)}$ are constant along the six different pieces.

We now first construct the solution $\Phi_{\alpha}^{(Bes)}$ of a model RH problem with the
same constant jumps on six infinite rays in an auxiliary $\zeta$-plane, and then we put
\begin{equation} \label{defPalphahat0}
    \widehat{P}_{\alpha}^{(0)}(z) = E_{\alpha}(z) \Phi_{\alpha}^{(Bes)}(t f(z)),
\end{equation}
where $f$ is given by
\begin{equation} f(z) = \frac{2}{3} - \frac{2i}{3}(z-1)^{3/2},
    \quad 0 < \arg(z-1) < 2 \pi,
\end{equation}
 and $E_{\alpha}(z)$ is an analytic prefactor that will be chosen later.
The function $f(z)$ is analytic in $U^{(0)}$. Moreover, $f(z) = z + \cdots$
and therefore it defines a conformal map in the neighborhood $U^{(0)}$.
After performing a slight contour deformation, we may and do assume that the
six contours $\Sigma_1$, $\Sigma_2$, $\Sigma_3$, $\Sigma_4$,
$\Sigma_u$ and $\Sigma_d$ are mapped by $f$ into six rays.
The exact relation between the functions $f(z)$ and $g(z)$ is given by the formula
\begin{equation}
    f(z) =
     \frac{2}{3} - i \left\{ \begin{array}{ll}
    g(z),& \Im z >0,\\[10pt]
   -g(z),& \Im z < 0.
\end{array} \right.
\end{equation}

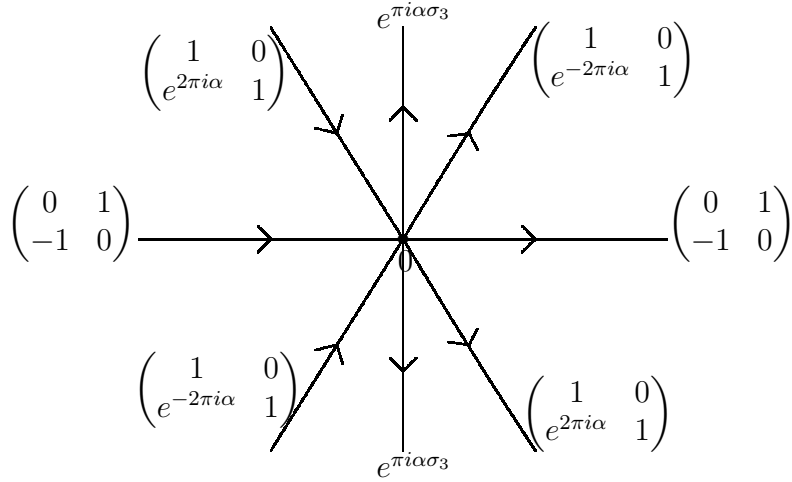
\begin{figure}[th]
\centering
\unitlength 1pt
\linethickness{0.5pt}
\begin{picture}(200,160)(-100,-80)
   \put(0,0){\line(1,0){100}}
   \put(0,0){\line(-1,0){100}}
   \qbezier(0,0)(-25,40)(-50,80)
   \qbezier(0,0)(-25,-40)(-50,-80)
   \qbezier(0,0)(30,50)(50,80)
   \qbezier(0,0)(30,-50)(50,-80)
   \qbezier(0,-80)(0,0)(0,80)
   \qbezier(50,0)(49,1)(45,5)   \qbezier(50,0)(49,-1)(45,-5)
   \qbezier(-50,0)(-51,1)(-55,5)   \qbezier(-50,0)(-51,-1)(-55,-5)
   \qbezier(-25,40)(-24,43.5)(-23,47)   \qbezier(-25,40)(-29,41)(-33,42)
   \qbezier(-25,-40)(-24,-43.5)(-23,-47)   \qbezier(-25,-40)(-29,-41)(-33,-42)
   \qbezier(25,40)(21,39.5)(17,39)   \qbezier(25,40)(26.5,37)(28,34)
   \qbezier(25,-40)(21,-39.5)(17,-39)   \qbezier(25,-40)(26.5,-37)(28,-34)
   \qbezier(0,50)(1,49)(5,45)   \qbezier(0,50)(-1,49)(-5,45)
   \qbezier(0,-50)(1,-49)(5,-45)   \qbezier(0,-50)(-1,-49)(-5,-45)
   \put(-2,-12){$0$}
   \put(0,0){\circle*{4}}
   \put(100,3){$\begin{pmatrix} 0& 1  \\ -1 & 0\end{pmatrix}$}
   \put(-100,60){$\begin{pmatrix} 1 & 0 \\  e^{2\pi i\alpha} & 1 \end{pmatrix}$}
   \put(-150,3){$\begin{pmatrix} 0& 1  \\ -1 & 0\end{pmatrix}$}
   \put(-102,-60){$\begin{pmatrix} 1 & 0 \\ e^{-2\pi i\alpha} & 1 \end{pmatrix}$}
   \put(47,65){$\begin{pmatrix} 1 & 0 \\  e^{-2\pi i\alpha} & 1 \end{pmatrix}$}
   \put(45,-69){$\begin{pmatrix} 1 & 0 \\ e^{2\pi i\alpha} & 1 \end{pmatrix}$}
   \put(-10,80){$e^{\pi i \alpha \sigma_3}$}
   \put(-10,-90){$e^{\pi i \alpha \sigma_3}$}
   \end{picture}
   \caption{Contours and jump matrices for the RH problem for $\widetilde{\Phi}_{\alpha}^{(Bes)}$.}\label{figure8}
\end{figure}

We construct $\Phi_{\alpha}^{(Bes)}$ by relating it to a model RH problem
constructed by Vanlessen \cite{Vanl} and also used in \cite{KV},
whose solution we denote here by $\widetilde{\Phi}_{\alpha}^{(Bes)}$.
\paragraph{Riemann-Hilbert problem for $\widetilde{\Phi}_{\alpha}^{(Bes)}$}
\begin{enumerate}
\item[\rm (a)] $\widetilde{\Phi}^{(Bes)}_{\alpha} : \mathbb{C} \setminus \Gamma \to
    \mathbb C^{2\times 2}$ is analytic, where $\Gamma$ is the union of the eight
    half rays shown in Figure~\ref{figure8}, namely
    \[ \Gamma = \{ \zeta \in \mathbb C \mid \arg \zeta \in \{0, \pm \pi/3, \pm \pi/2, \pm 2\pi/3, \pi\} \}. \]
\item[\rm (b)] $\widetilde{\Phi}_{\alpha,+}^{(Bes)} = \widetilde{\Phi}_{\alpha,-}^{(Bes)}
    v_{\widetilde{\Phi}^{(Bes)}_{\alpha}}$ on $\Gamma$, where the constant
    jump matrices $v_{\widetilde{\Phi}^{(Bes)}_{\alpha}}$ are indicated
    in Figure~\ref{figure8}.
\item[\rm (d)]
$\widetilde{\Phi}_{\alpha}^{(Bes)}(\zeta) = O\begin{pmatrix} |\zeta|^{\alpha} & |\zeta|^{\alpha} \\
    |\zeta|^{\alpha} & |\zeta|^{\alpha}
\end{pmatrix}$ as $\zeta \to 0$, if $-1/2 < \alpha < 0$; and

    $\widetilde{\Phi}_{\alpha}^{(Bes)}(\zeta) = \left\{ \begin{array}{ll}
    O\begin{pmatrix} |\zeta|^{\alpha} & |\zeta|^{-\alpha} \\
    |\zeta|^{\alpha} & |\zeta|^{-\alpha} \end{pmatrix}
    & \text{as $\zeta \to 0$, $\zeta \in \mathbb C \setminus \Gamma$ with $\pi/3 < |\arg \zeta| < 2 \pi/3$}, \\[10pt]
    O\begin{pmatrix} |\zeta|^{-\alpha} & |\zeta|^{-\alpha} \\
    |\zeta|^{-\alpha} & |\zeta|^{-\alpha} \end{pmatrix}
    & \text{as $\zeta \to 0$, $\zeta \in \mathbb C \setminus \Gamma $ with $\zeta$ elsewhere},
    \end{array} \right.$
    if $\alpha \geq 0$.
\end{enumerate}

We did not specify the asymptotic condition. A solution of the RH problem for
$\widetilde{\Phi}_{\alpha}^{(Bes)}$ was given in terms of Bessel functions
of orders $\alpha \pm \frac{1}{2}$.
There is a different expression in each sector. Let it suffice here to mention
the solution in the sector $\pi/2 < \arg \zeta < 2\pi/3$. We give it in
a form that is different from the one in \cite{Vanl}, where the modified
Bessel functions $I_{\alpha\pm\frac{1}{2}}$ and $K_{\alpha \pm \frac{1}{2}}$
are used.
We state it here in terms of the usual Bessel functions $J_{\alpha \pm \frac{1}{2}}$
and the Hankel function of first kind $H_{\alpha \pm \frac{1}{2}}^{(1)}$ as follows
\begin{multline} \label{tildePhi}
    \widetilde{\Phi}_{\alpha}^{(Bes)} (\zeta)=
    \sqrt{\pi} e^{-\frac{1}{4}\pi i} \zeta^{1/2}
    \begin{pmatrix} J_{\alpha+\frac{1}{2}}(\zeta) & \frac{1}{2} H_{\alpha+\frac{1}{2}}^{(1)}(\zeta) \\
    J_{\alpha-\frac{1}{2}}(\zeta) & \frac{1}{2} H_{\alpha-\frac{1}{2}}^{(1)}(\zeta)
    \end{pmatrix}, \qquad
      \pi/2 < \arg \zeta < 2\pi/3.
        \end{multline}
The asymptotics as $\zeta \to \infty$ in this sector is
\begin{multline}
    \widetilde{\Phi}_{\alpha}^{(Bes)}(\zeta)
        = \frac{1}{\sqrt{2}} \left[\begin{pmatrix} 1  & -i \\ -i & 1 \end{pmatrix}
            + \frac{\alpha}{2\zeta}
            \begin{pmatrix} - i (\alpha + 1) & \alpha + 1 \\
            - (\alpha - 1) & i (\alpha - 1) \end{pmatrix} + O\left(\frac{1}{\zeta^2}\right) \right] \\
            \times
        e^{\frac{1}{4} \pi i \sigma_3}  e^{\frac{1}{2}\pi i \alpha \sigma_3} e^{-i \zeta \sigma_3},
        \qquad \pi/2 < \arg \zeta < 2 \pi/3.
        \end{multline}

Then we define
\begin{multline} \label{defPhiBes}
    \Phi_{\alpha}^{(Bes)}(\zeta) =
    e^{-(\frac{1}{2} \pi i \alpha + \frac{1}{4} \pi i) \sigma_3}
    \frac{1}{\sqrt{2}} \begin{pmatrix} 1 & i \\ i & 1 \end{pmatrix}
    \widetilde{\Phi}_{\alpha}^{(Bes)}(\zeta) \\
    \times
    \left\{ \begin{array}{ll}
     e^{\pi i \alpha \sigma_3},
        & \text{ if $\Re \zeta > 0$, $\Im \zeta > 0$}, \\
    I & \text{ if $\Re \zeta < 0$}, \\
     e^{-\pi i \alpha \sigma_3},
        & \text{ if $\Re \zeta > 0$, $\Im \zeta < 0$}.
    \end{array} \right.
\end{multline}
It is then easy to check that $\Phi_{\alpha}^{(Bes)}$ is analytic across $i \mathbb R$
and has the jump matrices indicated in Figure~\ref{figure7}, but of course on
the contour $\Gamma \setminus i \mathbb R$. The behavior at $0$ remains unaffected
by the above transformation, while the constant prefactor $
    e^{-(\frac{1}{2} \pi i \alpha + \frac{1}{4} \pi i) \sigma_3}
    \frac{1}{\sqrt{2}} \begin{pmatrix} 1 & i \\ i & 1 \end{pmatrix}$ is chosen so that we obtain the precise
    asymptotics as $\zeta \to \infty$ given in item (c) below.

\paragraph{Riemann-Hilbert problem for $\Phi^{(Bes)}_{\alpha}$}
\begin{enumerate}
\item[\rm (a)] $\Phi^{(Bes)}_{\alpha} : \mathbb{C} \setminus (\Gamma \setminus i \mathbb R) \to
    \mathbb C^{2\times 2}$ is analytic.
\item[\rm (b)] $\Phi_{\alpha,+}^{(Bes)} = \Phi_{\alpha,-}^{(Bes)}
    v_{\Phi^{(Bes)}_{\alpha}}$ on $\Gamma \setminus i \mathbb R$, where the constant
    jump matrices $v_{{\Phi}^{(Bes)}_{\alpha}}$ are indicated
    in Figure~\ref{figure7}.
\item[\rm (c)] $\Phi^{(Bes)}_{\alpha}(\zeta) =
\left(I + O\left(\frac{1}{\zeta}\right)\right)
    e^{-i \zeta \sigma_3}$ as $\zeta \to \infty$ with $\Im \zeta > 0$, and

     $\Phi^{(Bes)}_{\alpha}(\zeta) \begin{pmatrix} 0 & 1 \\ -1 & 0 \end{pmatrix} =
    \left(I + O\left(\frac{1}{\zeta}\right)\right)
    e^{-i \zeta \sigma_3}$ as $\zeta \to \infty$ with $\Im \zeta < 0$.

\item[\rm (d)]
$\Phi_{\alpha}^{(Bes)}(\zeta) = O\begin{pmatrix} |\zeta|^{\alpha} & |\zeta|^{\alpha} \\
    |\zeta|^{\alpha} & |\zeta|^{\alpha}
\end{pmatrix}$ as $\zeta \to 0$, if $-1/2 < \alpha < 0$; and

    $\Phi_{\alpha}^{(Bes)}(\zeta) = \left\{ \begin{array}{ll}
    O\begin{pmatrix} |\zeta|^{\alpha} & |\zeta|^{-\alpha} \\
    |\zeta|^{\alpha} & |\zeta|^{-\alpha} \end{pmatrix}
    & \text{as $\zeta \to 0$, $\zeta \in \mathbb C \setminus \Gamma$ with $\pi/3 < |\arg \zeta| < 2 \pi/3$}, \\[10pt]
    O\begin{pmatrix} |\zeta|^{-\alpha} & |\zeta|^{-\alpha} \\
    |\zeta|^{-\alpha} & |\zeta|^{-\alpha} \end{pmatrix}
    & \text{as $\zeta \to 0$, $\zeta \in \mathbb C \setminus \Gamma $ with $\zeta$ elsewhere},
    \end{array} \right.$
    if $\alpha \geq 0$.
\end{enumerate}

As noted before, we now put
\[ P_{\alpha}^{(0)}(z) = E_{\alpha}(z) \Phi_{\alpha}^{(Bes)}(tf(z)) e^{tg(z) \sigma_3}, \]
where $E_{\alpha}$ is still to be determined.
Then, for fixed $z \in \partial U^{(0)}$ we have as $t \to \infty$,
\begin{align} \nonumber
    P_{\alpha}^{(0)}(z) & = E_{\alpha}(z) \Phi_{\alpha}^{(Bes)}(tf(z)) e^{tg(z)\sigma_3} \\
    \nonumber
    & = \left\{ \begin{array}{ll}
        E_{\alpha}(z) \left( I + O(1/t)\right) e^{-2 i t \sigma_3/3}, & \Im z > 0, \\
        E_{\alpha}(z) \begin{pmatrix} 0 & -1 \\ 1 & 0 \end{pmatrix}
            \left( I + O(1/t) \right) e^{2i t\sigma_3/3}, & \Im z < 0,
            \end{array} \right. \\
    \label{P0match1}
    & = \left\{ \begin{array}{ll}
        E_{\alpha}(z) e^{-2i t \sigma_3/3} \left( I + O(1/t)\right), & \Im z > 0, \\
        E_{\alpha}(z) e^{-2i t \sigma_3/3} \begin{pmatrix} 0 & -1 \\ 1 & 0 \end{pmatrix}
            \left( I + O(1/t) \right), & \Im z < 0.
            \end{array} \right.
     \end{align}
To match this with $P^{(\infty)}(z)$ for $z \in \partial U^{(0)}$ we choose
\begin{align} \label{defEalpha}
    E_{\alpha}(z) = \left\{ \begin{array}{ll}
    P^{(\infty)}(z) e^{2it\sigma_3/3}, & \quad \Im z > 0, \\
    P^{(\infty)}(z) \begin{pmatrix} 0 & 1 \\ -1 & 0 \end{pmatrix}
        e^{2it\sigma_3/3}, & \quad \Im z < 0,
        \end{array} \right.
        \end{align}
which is indeed analytic in $U^{(0)}$. This completes the construction
of the local parametrix $P_{\alpha}^{(0)}$.

\subsection{Fourth transformation $C_{\alpha} \mapsto D_{\alpha}$}
\begin{figure}[t]
\centering
\unitlength 1pt
\linethickness{0.5pt}
\begin{picture}(200,160)(-60,-80)
   \put(75,0){\line(1,0){135}}
   \qbezier(-8.5,13.6)(-25,40)(-50,80)
   \qbezier(-8.5,-13.6)(-25,-40)(-50,-80)
   \qbezier(8.5,13.6)(30,37)(51.5,13.6)
   \qbezier(8.5,-13.6)(30,-37)(51.5,-13.6)
   \qbezier(120,0)(119,1)(115,5)   \qbezier(120,0)(119,-1)(115,-5)
   \qbezier(-25,40)(-24,43.5)(-23,47)   \qbezier(-25,40)(-29,41)(-33,42)
   \qbezier(-25,-40)(-24,-43.5)(-23,-47)   \qbezier(-25,-40)(-29,-41)(-33,-42)
   \qbezier(30,25)(29,26)(25,30) \qbezier(30,25)(29,24)(25,20)
   \qbezier(30,-25)(29,-26)(25,-30) \qbezier(30,-25)(29,-24)(25,-20)

   \qbezier(-15,0)(-16,1)(-20,5)   \qbezier(-15,0)(-14,1)(-10,5)
   \qbezier(45,0)(44,1)(40,5)   \qbezier(45,0)(46,1)(50,5)
   \put(-2,-10){$0$}
   \put(0,0){\circle*{4}}
   \put(0,0){\circle{30}}
   \put(60,-10){$1$}
   \put(60,0){\circle*{4}}
   \put(60,0){\circle{30}}
   \put(-25,50){$\Sigma_{2}$}
   \put(-25,-55){$\Sigma_{4}$}
   \put(40,25){$\Sigma_u$}
   \put(45,-30){$\Sigma_d$}
   \end{picture}
   \caption{Contour $\Sigma_D$ in the RH problem for $D_{\alpha}$.}\label{figure9}
\end{figure}
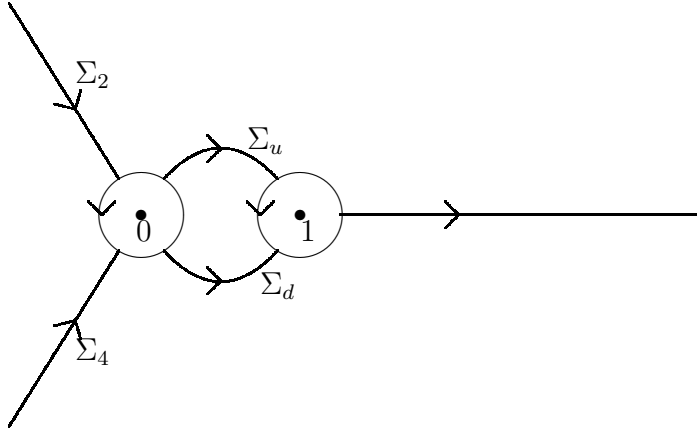

In the final transformation we put
\begin{equation} \label{Ddef2}
    D_{\alpha}(z) = \left\{ \begin{array}{ll}
        C_{\alpha}(z) \left(P^{(1)}(z)\right)^{-1}, & \text{for }
        z \in U^{(1)} \setminus \Sigma_C,\\[10pt]
        C_{\alpha}(z) \left(P^{(0)}_{\alpha}(z)\right)^{-1}, & \text{for }
        z \in U^{(0)} \setminus \Sigma_C,\\[10pt]
        C_{\alpha}(z) \left(P^{(\infty)}(z)\right)^{-1}, & \text{for }
        z \in \mathbb{C} \setminus (\overline{U^{(1)}\cup U^{(0)}  \cup (-\infty, 1)}).
        \end{array} \right.
\end{equation}
By construction, the only jumps that remain for $D_{\alpha}$ are
across the circles $\partial U^{(1)}$ and $\partial U^{(0)}$ and the parts of
the arcs $\Sigma_{u}$ and $\Sigma_{d }$ and the
rays $[1, \infty)$, $\Sigma_{2}$, and  $\Sigma_{4}$ which lie outside
of the neighborhoods $U^{(1)}$ and $U^{(0)}$. We shall denote this
remaining contour as $\Sigma_{D}$; it is depicted in
Figure~\ref{figure9}. The Riemann-Hilbert problem for $D_{\alpha}$ is set
on this contour.

\paragraph{Riemann-Hilbert problem for $D_{\alpha}$}
\begin{itemize}
\item[(a)] $D_{\alpha} : \mathbb{C} \setminus \Sigma_{D} \to \mathbb C^{2\times 2}$ is analytic.
\item[(b)] $D_{\alpha,+}(z) = D_{\alpha,-}(z)\,v_{D_{\alpha}}(z)$ for $z \in \Sigma_{D}$, where
\begin{equation} \label{vRdef2}
v_{D_{\alpha}} = \left\{
    \begin{array}{ll}
        P^{(\infty)}\,(P^{(1)})^{-1},&\text{ on } \partial U^{(1)},\\
        P^{(\infty)}\,(P^{(0)})^{-1},&\text{ on } \partial U^{(0)},\\
        P^{(\infty)}\,v_{C_{\alpha}}\,(P^{(\infty)})^{-1},&\text{ on }
        \Sigma_{D} \setminus (\partial U^{(1)}\cup\partial U^{(0)}) .
\end{array}
\right.
\end{equation}
\item[(c)] $D_{\alpha}(z) = I + O(1/z)$ as $z \to \infty$.
\end{itemize}
Due to the matching conditions of the Riemann-Hilbert
problems for $P^{(1)}$ and $P^{(0)}$, we have that
\begin{equation}\label{jump1}
    v_{D_{\alpha}}(z) = I +  O\left(\frac{1}{t}\right),
\end{equation}
uniformly on the circles   $\partial U^{(1)}$ and $\partial U^{(0)}$.
Simultaneously,
\begin{equation}\label{jump2}
    v_{D_{\alpha}}(z) = I +  O\left(e^{-ct|z|}\right), \quad c >0,
\end{equation}
uniformly on $ \Sigma_{D} \setminus (\partial U^{(1)}\cup\partial U^{(0)})$.
Hence, as before,
\begin{equation}\label{Das}
D_{\alpha}(z) = I + O\Big(\frac{1}{t(1 + |z|)}\Big), \quad \text{ as } t \to +\infty,
\end{equation}
uniformly for $z \in \mathbb{C}\setminus \Sigma_D$.

\subsection{Conclusion of the proof of \eqref{uminusinfinity}}

The main remaining step in the proof of \eqref{uminusinfinity}
is to express $u_{\alpha}$ in terms of $D_{\alpha}$. The result of the
calculations is contained in the next lemma.
\begin{lemma} \label{lemualphainD}
For every $s < 0$, we have
\begin{equation} \label{usolution5D}
    u_{\alpha}(s) = \frac{\alpha}{\sqrt{-s}}
    \left[ D_{\alpha}(0) \begin{pmatrix} i \sin \theta & \cos \theta \\
        -\cos \theta & - i \sin \theta \end{pmatrix} D_{\alpha}^{-1}(0) \right]_{12}.
        \end{equation}
where
\begin{equation} \label{deftheta}
    \theta = 4 t/3 - \pi \alpha.
    \end{equation}
and $t = (-s)^{3/2}$ as before.
\end{lemma}

Theorem \ref{theorem2} follows immediately from the lemma and \eqref{Das}.
Indeed from \eqref{Das} we find that
$D_{\alpha}(0) = I + O((-s)^{-3/2})$ as $s \to -\infty$.
Hence by \eqref{usolution5D} and \eqref{deftheta}
\begin{align*} u_{\alpha}(s) & = \frac{\alpha}{\sqrt{-s}}
        \left( \cos \theta + O((-s)^{-3/2}) \right)  = \frac{\alpha}{\sqrt{-s}}
            \cos\left( \frac{4}{3} (-s)^{3/2}  - \pi \alpha\right) + O(s^{-2})
            \end{align*}
            as $s \to -\infty$, which is \eqref{uminusinfinity}.

So it remains to prove Lemma \ref{lemualphainD}.
\begin{proof}
Take $z \in U^{(0)}$ with $\Im z >0 $ and outside of the lense.
Then we have by \eqref{Ddef2}, \eqref{defPalpha0}, and \eqref{defPalphahat0}, 
that
\begin{align} \nonumber
    C_{\alpha}(z) & = D_{\alpha}(z) P_{\alpha}^{(0)}(z) \\
        & = D_{\alpha}(z) E_{\alpha}(z)
        \Phi_{\alpha}^{(Bes)}(t f(z))
                e^{t g(z) \sigma_3}. \label{fromCtoD}
                \end{align}
In the computation of the limit of $z\left[ \left( \frac{d}{dz} C_{\alpha}(z) \right) C_{\alpha}^{-1}(z) \right]_{12}$
as $z \to 0$, the only term that will contribute is the one
we get by taking the derivative of $\Phi_{\alpha}^{(Bes)}(tf(z))$.
This easily follows from the fact that $D_{\alpha}$ and $E_{\alpha}$ are analytic at $0$.

Since, by \eqref{defEalpha} and \eqref{Pinfsol}, we have
\begin{align} \label{Ealphaat0}
    E_{\alpha}(0) = e^{-i \pi \sigma_3/4}
    \frac{1}{\sqrt{2}} \begin{pmatrix} 1 & i \\ i & 1 \end{pmatrix} e^{2it \sigma_3/3}
\end{align}
    it then follows from
\eqref{usolution4C} and \eqref{fromCtoD} that
\begin{align} \nonumber
    u_{\alpha}(s) = & \frac{i}{\sqrt{-s}}
    \left[ D_{\alpha}(0) e^{-i \pi \sigma_3/4}
        \frac{1}{\sqrt{2}} \begin{pmatrix} 1 & i \\ i & 1 \end{pmatrix}
        e^{2i t \sigma_3/3} \right. \\
           & \nonumber \qquad \times \left(    \lim_{z \to 0}
          z  \left(\frac{d}{dz} \Phi_{\alpha}^{(Bes)}(tf(z))\right)
          \left( \Phi_{\alpha}^{(Bes)}(tf(z))\right)^{-1} \right) \\
          & \qquad \qquad \times \left.
        e^{-2i t \sigma_3/3}
        \frac{1}{\sqrt{2}} \begin{pmatrix} 1 & -i \\ -i & 1 \end{pmatrix}
        e^{i \pi \sigma_3/4} D_{\alpha}^{-1}(0) \right]_{12}.
        \label{usolution5Da}
        \end{align}

Putting $\zeta = tf(z)$ we get
\[           \frac{d}{dz} \Phi_{\alpha}^{(Bes)}(tf(z)) =
           t f'(z) \frac{d}{d\zeta} \Phi_{\alpha}^{(Bes)}(\zeta). \]
Noting that $f(0) = 0$ and $f'(0) \neq 0$, we find
\[ \frac{t f'(z) z}{\zeta} = \frac{z f'(z)}{f(z)} \to 1 \]
as $z \to 0$. Therefore
\begin{equation} \label{logBesselA}
    \lim_{z \to 0} z
          \left( \frac{d}{dz} \Phi_{\alpha}^{(Bes)}(tf(z))\right) \left( \Phi_{\alpha}^{(Bes)}(tf(z))\right)^{-1}
  = \lim_{\zeta \to 0} \zeta
           \left( \frac{d}{d\zeta} \Phi_{\alpha}^{(Bes)}(\zeta)\right) \left( \Phi_{\alpha}^{(Bes)}(\zeta)\right)^{-1}.
            \end{equation}

From the definition \eqref{defPhiBes} of $\Phi_{\alpha}^{(Bes)}$
we find that for $\Im \zeta > 0$,
\begin{multline}  \label{logBesselB}
   \left( \frac{d}{d\zeta} \Phi_{\alpha}^{(Bes)}(\zeta) \right) \left( \Phi_{\alpha}^{(Bes)}(\zeta)\right)^{-1} \\
    =
    e^{-( \pi i \alpha/2 + \pi i/4) \sigma_3}
        \frac{1}{\sqrt{2}} \begin{pmatrix} 1 & i \\ i & 1 \end{pmatrix}
    \left( \left( \frac{d}{d\zeta} \widetilde{\Phi}_{\alpha}^{(Bes)}(\zeta) \right)
    \left( \widetilde{\Phi}_{\alpha}^{(Bes)}(\zeta)\right)^{-1} \right) \\
    \times  \frac{1}{\sqrt{2}} \begin{pmatrix} 1 & -i \\ -i & 1 \end{pmatrix}
    e^{(\pi i \alpha/2 + \pi i/4) \sigma_3}.
    \end{multline}
Recall that $\widetilde{\Phi}_{\alpha}^{(Bes)}$ is built out of
Bessel functions. From the differential-difference relations
satisfied by the Bessel functions (see formula 9.1.27 in
\cite{AS})
\begin{align*}
    J_{\nu}'(\zeta) & = J_{\nu-1}(\zeta) - \frac{\nu}{\zeta} J_{\nu}(\zeta), \\
    J_{\nu-1}'(\zeta) & = - J_{\nu+1}(\zeta) + \frac{\nu}{\zeta} J_{\nu}(\zeta),
    \end{align*}
and similar ones for the Hankel functions, it easily follows from \eqref{tildePhi} that
\begin{align*}
    \frac{d}{d\zeta} \widetilde{\Phi}_{\alpha}^{(Bes)}(\zeta) =
    \begin{pmatrix} -\alpha/\zeta & 1 \\ -1 & \alpha/\zeta \end{pmatrix}
         \widetilde{\Phi}_{\alpha}^{(Bes)}(\zeta).
         \end{align*}
Thus
\begin{equation} \label{logBesselC} \lim_{\zeta \to 0} \zeta
    \left( \left( \frac{d}{d\zeta} \widetilde{\Phi}_{\alpha}^{(Bes)}(\zeta) \right)
     \left(\widetilde{\Phi}_{\alpha}^{(Bes)}(\zeta) \right)^{-1} \right)
     = - \alpha \sigma_3. \end{equation}

Combining \eqref{logBesselC} with \eqref{usolution5Da}, \eqref{logBesselA}, \eqref{logBesselB}
we get with $\theta = 4 t/3 - \pi \alpha$,
\begin{align*} u_{\alpha}(s) = & \frac{i}{\sqrt{-s}}
    \left[ D_{\alpha}(0) e^{-i \pi \sigma_3/4}
        \frac{1}{\sqrt{2}} \begin{pmatrix} 1 & i \\ i & 1 \end{pmatrix}
        e^{ \theta i \sigma_3/2 - \pi i \sigma_3/4}
        \frac{1}{\sqrt{2}} \begin{pmatrix} 1 & i \\ i & 1 \end{pmatrix} \right.  \\
        & \times \left.
        (-\alpha \sigma_3)
        \frac{1}{\sqrt{2}} \begin{pmatrix} 1 & -i \\ -i & 1 \end{pmatrix}
        e^{-\theta i \sigma_3/2 + \pi i \sigma_3/4}
        \frac{1}{\sqrt{2}} \begin{pmatrix} 1 & -i \\ -i & 1 \end{pmatrix}
        e^{i \pi \sigma_3/4} D_{\alpha}^{-1}(0) \right]_{12},
        \end{align*}
which after straightforward calculation reduces to \eqref{usolution5D}.
This completes the proof of the lemma.
\end{proof}

\section*{Acknowledgements}

Alexander Its was supported in part by NSF grant \#DMS-0401009.
Arno Kuijlaars is supported by FWO-Flanders project G.0455.04,
by K.U.~Leuven research grants OT/04/21 and OT/08/33, by the Belgian Interuniversity
Attraction Pole P06/02, by the  European Science Foundation Program MISGAM,
and by a grant from the Ministry of Science and Innovation of Spain,
project code MTM2008-06689-C02-01.
J\"{o}rgen \"{O}stensson is supported by K.U.~Leuven research
grant OT/04/21.

\appendix

\section{Appendix. Relation to the Painlev\'e II equation and the uniqueness question}
\label{sectionA}

We start with  reviewing the general facts concerning the
relation between the thirty fourth and the second Painlev\'e equations.

Let $q(s)$ be a solution of the second Painlev\'e equation with
parameter $\nu$,
\begin{equation}\label{painleve2}
q'' = 2q^3 + sq -\nu.
\end{equation}
Then, the function $u(s)$ defined by the formulae
\begin{equation}\label{p2top34}
 u(s) = 2^{-1/3} U(-2^{1/3} s), \qquad  U(s) = q^2(s) + q'(s) + \frac{s}{2}.
\end{equation}
satisfies the thirty fourth equation \eqref{painleve34} with the parameter
\[
\alpha = \frac{\nu}{2} -\frac{1}{4}
\]
(see \cite{Ince};
see also \cite{BBIK} and \cite{IKO1}).
The inverse transformation is given by the formulae
\begin{equation}\label{p34top2}
 q(s) = -2^{-1/3} Q(-2^{-1/3} s), \qquad  Q(s) = \frac{u' -2\alpha}{2u}.
\end{equation}
Moreover, the equations (see \cite{BBIK}, \cite{KH}),
\begin{equation}
\label{PsiFNtoPsialpha}
\Psi_{\alpha}(z;s) = \begin{pmatrix} 1 & 0 \\ \eta(s) & 1 \end{pmatrix}
    z^{-\sigma_3/4} \frac{1}{\sqrt{2}} \begin{pmatrix}
    1 & i \\ i & 1 \end{pmatrix} e^{\pi i \sigma_3/4}
    \Psi_{2\alpha + 1/2}^{FN}(w;-2^{1/3}s) e^{-\pi i \sigma_3/4},
\end{equation}
\[
\eta(s) = -2^{1/3}\left([m_{FN}(s)]_{11} + [m_{FN}(s)]_{12}\right),
\]
where $w = e^{\pi i/2} 2^{-1/3} z^{1/2}$ with $\Im w > 0$,
establish the relation between the solution $\Psi_{\alpha}(z;s)$
of the general Painlev\'e XXXIV  RH problem formulated in
Section \ref{RHproblem34}  and the
solution  $\Psi_{\nu}^{FN}(w;s)$ of the
RH problem associated with  the Painlev\'e II equation \cite{FN}
with the parameter $\nu = 2\alpha + 1/2$.
In \eqref{PsiFNtoPsialpha},
$m_{FN}(s)$ denotes the first matrix coefficient
of the expansion  $\Psi_{\nu}^{FN}(w;s)$ at $w =\infty$,
\begin{equation} \label{PFNsiexpansion}
    \Psi^{FN}_{\nu}(w;s) =
    \left(I + \frac{m_{FN}(s)}{z} + O\left(\frac{1}{w^2}\right)\right)
e^{-i(\frac{4}{3} w^{3} + s w)\sigma_3}
\end{equation}
as $w \to \infty$.
We use here the Flaschka-Newell \cite{FN} form of the Painlev\'e II RH problem whose setting
we will now remind (for details see \cite[Chapter 5]{FIKN}).

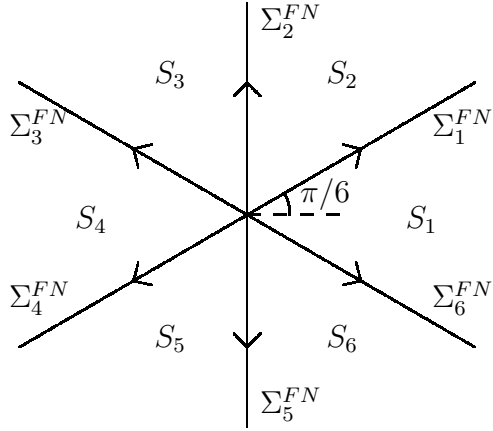
\begin{figure}[th]
\centering
\unitlength 1pt
\linethickness{0.5pt}
\begin{picture}(200,160)(-100,-80)
   \put(0,0){\line(1,0){5}}
   \put(10,0){\line(1,0){5}}
   \put(20,0){\line(1,0){5}}
   \put(30,0){\line(1,0){5}}
   \qbezier(0,0)(-43,25)(-86,50)
   \qbezier(0,0)(-43,-25)(-86,-50)
   \qbezier(0,0)(43,25)(86,50)
   \qbezier(0,0)(43,-25)(86,-50)
   \qbezier(0,-80)(0,0)(0,80)
   \qbezier(-43,25)(-43,25)(-41.2,18)   \qbezier(-43,25)(-43,25)(-36.2,26.6)
   \qbezier(-43,-25)(-43,-25)(-41.2,-18)   \qbezier(-43,-25)(-43,-25)(-36.2,-26.6)
   \qbezier(43,25)(43,25)(41.2,18)   \qbezier(43,25)(43,25)(36.2,26.6)  
   \qbezier(43,-25)(43,-25)(41.2,-18)   \qbezier(43,-25)(43,-25)(36.2,-26.6)
   \qbezier(0,50)(1,49)(5,45)   \qbezier(0,50)(-1,49)(-5,45)
   \qbezier(0,-50)(1,-49)(5,-45)   \qbezier(0,-50)(-1,-49)(-5,-45)
    \put(70,30){$\Sigma_1^{FN}$}
    \put(5,70){$\Sigma_2^{FN}$}
    \put(-90,30){$\Sigma_3^{FN}$}
    \put(-90,-35){$\Sigma_4^{FN}$}
    \put(5,-75){$\Sigma_5^{FN}$}
    \put(70,-35){$\Sigma_6^{FN}$}
    \put(60,-5){$S_1$}
    \put(30,50){$S_2$}
    \put(-35,50){$S_3$}
    \put(-65,-5){$S_4$}
    \put(-35,-50){$S_5$}
    \put(30,-50){$S_6$}
    \put(20,5){$\pi/6$}
    \qbezier(16,0)(16.5,4.4)(14,8)
   \end{picture}
   \caption{Contour for the RH problem for $\Psi_{\nu}^{FN}$.} \label{figure10}
\end{figure}

The general Painlev\'e II  RH problem involves three complex constants $a_1$, $a_2$, $a_3$ satisfying
(cf. \eqref{b_s})
\begin{equation} \label{a1a2a3variety}
    a_1 + a_2 + a_3 + a_1 a_2 a_3 = - 2i \sin \nu \pi,
    \end{equation}
and certain connection matrices $E^{FN}_{j}$.
 Let $S_j = \{w \in \mathbb{C} \mid \frac{2j-3}{6}\pi < \arg w <
\frac{2j-1}{6}\pi \}$ for $j = 1,\ldots,6$, and let $\Sigma^{FN} =
\mathbb C \setminus \bigcup_j S_j$. Then $\Sigma^{FN}$ consists of
six rays $\Sigma_j^{FN} = \{w \in \mathbb{C} \mid \arg w =
\frac{2j-1}{6}\pi \} $ for $j=1,\ldots,6$, all chosen oriented
towards infinity as in Figure~\ref{figure10}.

The RH problem is the following.
\begin{enumerate}
\item[\rm (a)] $\Psi_{\nu}^{FN} : \mathbb{C} \setminus \Sigma^{FN}  \to
    \mathbb C^{2\times 2}$ is analytic,
\item[\rm (b)] $\Psi_{\nu,+}^{FN} = \Psi_{\nu,-}^{FN}
    \begin{pmatrix} 1 & 0 \\ a_1 & 1 \end{pmatrix}$
    on $\Sigma_1^{FN}$,

    $\Psi_{\nu,+}^{FN} = \Psi_{\nu,-}^{FN}
    \begin{pmatrix} 1 & a_2 \\ 0 & 1 \end{pmatrix}$
    on $\Sigma_2^{FN}$,

    $\Psi_{\nu,+}^{FN} = \Psi_{\nu,-}^{FN}
    \begin{pmatrix} 1 & 0 \\ a_3 & 1 \end{pmatrix}$
    on $\Sigma_3^{FN}$,

    $\Psi_{\nu,+}^{FN} = \Psi_{\nu,-}^{FN}
    \begin{pmatrix} 1 & a_1 \\ 0 & 1 \end{pmatrix}$
    on $\Sigma_4^{FN}$,

    $\Psi_{\nu,+}^{FN} = \Psi_{\nu,-}^{FN}
    \begin{pmatrix} 1 & 0 \\ a_2 & 1 \end{pmatrix}$
    on $\Sigma_5^{FN}$,

    $\Psi_{\nu,+}^{FN} = \Psi_{\nu,-}^{FN}
    \begin{pmatrix} 1 & a_3 \\ 0 & 1 \end{pmatrix}$
    on $\Sigma_6^{FN}$.
\item[\rm (c)] $\Psi_{\nu}^{FN}(w) = (I + O(1/w)) e^{-i(\frac{4}{3} w^3 + s w) \sigma_3}$
    as $w \to \infty$.
\item[\rm (d)] If $\nu - \frac{1}{2} \not\in \mathbb N_0$, then
\begin{equation} \label{connection0}
    \Psi_{\nu}^{FN}(w) = B(w)
    \begin{pmatrix} w^{\nu} & 0 \\ 0 & w^{-\nu} \end{pmatrix} E^{FN}_j,
    \quad \text{ for } w \in S_j,
\end{equation}
where $B$ is analytic.
If $\nu \in \frac{1}{2} + \mathbb N_0$, then
there exists a constant $\kappa$ such that
\begin{equation} \label{connection1}
    \Psi_{\nu}^{FN}(w) = B(w)
    \begin{pmatrix} w^{\nu} & \kappa w^{\nu}\log w \\ 0 & w^{-\nu} \end{pmatrix} E^{FN}_j,
    \quad \text{ for } w \in S_j,
\end{equation}
where $B$ is analytic.
\end{enumerate}
Except for the special case,
\begin{equation} \label{specialcase}
    \nu = \frac{1}{2} + n, \quad a_1 = a_2 = a_3 = i (-1)^{n+1},
\quad n \in \mathbb N, \end{equation}
when the solution of the RH problem is  given in fact in terms of the Airy functions,
the connection matrices $E^{FN}_{j}$ are determined
(up to inessential left diagonal or upper triangular factors) by $\nu$ and
the Stokes multipliers $a_{j}$. In the special case the solution
is parametrized by the one non-trivial entry of the connection
matrix $E^{FN}_{1}$. We refer to  \cite[Chapters 5 and 11]{FIKN} for more details on
the setting and the analysis of the Painlev\'e II RH problem (see also
our paper \cite{IKO1}, where we review these results in the notations
we use here).

Equation \eqref{PsiFNtoPsialpha} implies
the following relation between the Painlev\'e XXXIV
and Painlev\'e II Stokes
parameters (\cite{BBIK}, \cite{KH}; see also \cite{IKO1}),
\begin{equation}\label{atob}
b_{1} = ia_{2},\quad b_{2} = ia_{3},\quad b_{4} = ia_{1},
\end{equation}
where we use the $b_j$ as in Section
\ref{RHproblem34}.
Therefore, taking into account \eqref{p34spec},
we conclude that  the second Painlev\'e function which is related to
the special solution $u_{\alpha}(s)$ of the thirty fourth Painlev\'e
equation studied in this paper corresponds to the choice
of the Stokes multipliers,
\begin{equation}\label{specP2}
    a_{1} = e^{-\nu\pi i}, \quad a_{2} =-i,\quad a_{3} = -e^{\nu\pi i}.
\end{equation}
This is different from the choice,
\begin{equation}\label{HM}
    a_{1} =  e^{-\nu\pi i}, \quad a_{2} = 0,\quad
    a_{3} = -e^{ \nu\pi i},
\end{equation}
corresponding, as it is shown in \cite[Chapter 11]{FIKN}
\footnote{The  RH problem which is used in \cite[Chapter 11]{FIKN} differs
by a simple gauge transformation from the Flaschka-Newell RH problem.
Indeed one has that $\Psi^{FN} = e^{-i\frac{\pi}{4}\sigma_{3}}\Psi e^{i\frac{\pi}{4}\sigma_{3}}$,
where $\Psi(z)$ is the solution of the RH problem from \cite[Chapter 11]{FIKN}.
This in turn implies that the Flaschka-Newell monodromy parameters $a_{j}$
are related to the monodromy parameters $s_{j}$ from \cite[Chapter 11]{FIKN},
via the equations, $a_1 = is_{1}$, $a_2 = -is_{2}$, and $a_3 = is_{3}$.},
to the generalized ($\nu \neq 0$) Hastings-McLeod
solution of the second Painlev\'e equation, i.e., the solution which
is characterized by the following asymptotic conditions,
\begin{equation}\label{HMminus}
q_{HM}(s) =  \sqrt{-\frac{s}{2}} + O(s^{-1}),
\qquad \text{as } s \to -\infty,
\end{equation}
and
\begin{equation}\label{HMplus}
q_{HM}(s) = \frac{\nu}{s} +  O(s^{-4}),
\qquad \text{as } s \to +\infty,
\end{equation}
see also \cite{CKV}.

It follows from \eqref{specP2}, however, that both solutions - the Hastings-McLeod
solution and the one corresponding to  $u_{\alpha}(s)$, belong
to the same one-parameter family of solutions of the Painlev\'e II
equation which is characterizes by the following choice
of the Stokes multipliers:
\begin{equation}\label{specP2a2free}
a_{1} =  e^{- \nu\pi i}, \quad
a_{3} = -e^{ \nu\pi i},
\end{equation}
and the Stokes multiplier $a_{2}$ is a free parameter of the family.
It is shown in \cite{Kap2} (see also \cite[Chapter 11]{FIKN})
that this family is exactly the classical  family of the so called {\it tronqu\'ee}
solutions, i.e., the solutions all of  which exhibit the same
behavior \eqref{HMminus} at $-\infty$.
In fact, for every {\it tronqu\'ee}  solution the
behavior \eqref{HMminus}
can be extended to a full asymptotic series and it holds in the
whole sector $2\pi/3 < \arg s < 4\pi/3$,
\begin{alignat}{2}
\nonumber
q_{tronq}(s) & \sim \sqrt{-\frac{s}{2}}\sum_{n=0}^{\infty}c_n(-s)^{-3n/2},
\quad c_{0} =1,\\
\label{tritrongP2}
& \qquad \text{ as } s \to -\infty, \quad \arg s \equiv \pi + \arg(-s)
\in \bigl(\tfrac{2\pi}{3}, \tfrac{4\pi}{3}\bigr),
\end{alignat}
where all the coefficients $c_{n}$ are uniquely determined by the
substitution into the Painlev\'e II equation (i.e., the series is
the {\it same} for every solution from the family
\footnote{We refer to \cite[Chapter 11]{FIKN} for more
on the asymptotic analysis of the {\it tronqu\'ee} solutions of the second
Painlev\'e equation. In particular, the reader can find there an alternative
parametrization (and its explicit relation to $a_{2}$)  of the
solutions via the coefficients of the oscillatory terms of the asymptotics
on the boundary rays.}).

The above made observation suggests that one can obtain the asymptotic
statement \eqref{uplusinfinity} of Theorem \ref{theorem2} directly
from \eqref{tritrongP2} using relation \eqref{p2top34}. Indeed,
the first two terms of  \eqref{tritrongP2} read
\begin{alignat}{2}
\nonumber
q_{tronq}(s) & =  \sqrt{-\frac{s}{2}} -\frac{\nu}{2s} + O(s^{-5/2}),\\
\label{tritrongP22}
& \qquad \text{ as } s \to -\infty, \quad \arg s \equiv \pi + \arg(-s)
\in \bigl(\tfrac{2\pi}{3}, \tfrac{4\pi}{3}\bigr).
\end{alignat}
Substituting this (differentiable !) asymptotics into \eqref{p2top34}
we indeed arrive at \eqref{uplusinfinity}.

The fact that the asymptotics \eqref{uplusinfinity} holds for a one-parameter
family of the solutions of the Painlev\'e XXXIV equation \eqref{painleve34}
can also be deduced from the direct analysis of the Painlev\'e XXXIV
RH problem. Denote
the one-parameter  family of solutions of Painlev\'e XXXIV corresponding to the
RH data,
\begin{equation}\label{btrongp34}
b_{2} = e^{2\alpha \pi i} ,\quad b_{4} = e^{-2\alpha \pi i},
\quad b_{1} = b \in  \mathbb C,
\end{equation}
as
\[
u^{(tronq)}_{\alpha}(s) \equiv u^{(tronq)}_{\alpha}(s|b).
\]
Note, that the cyclic relation
\eqref{b_s} is valid identically for $b_{1}$ if $b_{2}$ and $b_{4}$ are
as in \eqref{btrongp34}. It is not difficult to see that
{\it exactly} the same sequence of transformation as the one we
used in Section \ref{section2} in the analysis of the RH problem
in the case $s \to + \infty$ can be performed for any value of $b$.
In the final $D_{\alpha}$ RH problem the only difference is
in the jump matrix $v_{D_{\alpha}}$ on the segment of the horizontal
part of the jump contour depicted in Figure~\ref{figure3} which is to the right of $0$.
That is, instead of
\[
v_{D_{\alpha}} = P_{\alpha}^{(\infty)}(z)
 \begin{pmatrix} 1 & e^{-2t g(z)} \\ 0 & 1 \end{pmatrix}
 \left(P_{\alpha}^{(\infty)}(z)\right)^{-1} \equiv I +
 O\left(e^{-ct(|z| +1)}\right),
 \qquad \text{for} \quad z >0,
\]
 we now have,
\[
v_{D_{\alpha}} = P_{\alpha}^{(\infty)}(z)
 \begin{pmatrix} 1 & be^{-2t g(z)} \\ 0 & 1 \end{pmatrix}
 \left(P_{\alpha}^{(\infty)}(z)\right)^{-1} \equiv I +
 O\left(be^{-ct(|z| +1)}\right),
 \qquad \text{for} \quad z >0.
\]
In other words, the only difference in $v_{D_{\alpha}}$ is in
the  exponentially small error. Hence the estimates
\eqref{vD1} and \eqref{vD2} are valid for all $b$ and lead
to the same asymptotic behavior  \eqref{uplusinfinity}
of the solution $u^{(tronq)}_{\alpha}(s)$ for all $b$.
In fact, every solution from the family has the same
asymptotic series representation,
\begin{equation}\label{trongP34}
u^{(tronq)}_{\alpha}(s) \sim \frac{\alpha}{\sqrt{s}} + \sum_{n=1}^{\infty}d_n s^{-\frac{3n+1}{2}},
\quad \text{ as } s \to  + \infty,
\end{equation}
with the coefficients $d_{n}$ uniquely determined by the substitution
of the series into equation \eqref{painleve34} (see \cite{IKO1} for
the explicit recurrence relation for $d_{n}$).

The asymptotic behavior \eqref{HMplus} of the Hastings-McLeod
solution at $+\infty$ is also shared by another one-parameter
family of {\it tronqu\'ee}  solutions. The corresponding RH parametrization
is (see \cite{IKap2}; see also \cite[Chapter 11]{FIKN}),
\begin{equation}\label{trongplus}
a_{2} = 0, \quad a_{1} + a_{3} = -2i \sin\nu\pi.
\end{equation}
The Painlev\'e II function which corresponds to the thirty fourth
Painlev\'e function $u_{\alpha}(s)$ obviously does not belong to this
family and hence does not behave as  \eqref{HMplus} when $s \to + \infty$.
However, the leading term of its behavior as $s \to +\infty$ is
known (\cite{Kap1}; see also \cite[Chapter 10]{FIKN}).
Unfortunately, the leading term is not enough to derive the
corresponding asymptotics as $s \to -\infty$ of the Painlev\'e
XXXIV function $u_{\alpha}(s)$. Indeed, the leading asymptotics of
$q(s)$ as $s \to +\infty$ is of the form
\begin{equation}
\label{P2atplusinfinity}
q(s) \sim \sqrt{\frac{s}{2}} \cot \left(\frac{\sqrt{2}}{3}s^{3/2} + \chi \right),
\end{equation}
(the phase $\chi$ is known) and it cancels out in the right-hand side of equation
\eqref{p2top34}. Moreover, the solution $q(s)$, as it follows from
\eqref{P2atplusinfinity}, has poles on the positive real $s$-axis while the
function $u_{\alpha}(s)$ is smooth for all real $s$. This means that
the reduction to the Painlev\'e II equation is not the
best way to study the asymptotics of the Painlev\'e XXXIV function
$u_{\alpha}(s)$ as $s \to -\infty$. It is better to proceed via
the direct analysis of the Painlev\'e XXXIV RH problem for $\Psi_{\alpha}$,
as we did in Section \ref{section3} of this paper
\footnote{In principle, it is possible to use the inverse formula
\eqref{p34top2} and obtain the asymptotic of $u_{\alpha}(s)$
by integrating \eqref{P2atplusinfinity}. However, with this approach we
face the problem of  evaluation of the constant of integration
and, once again, one has to take special care of the poles
of the Painlev\'e II function.}.

Let us now consider the question of the uniqueness of the
solution $u_{\alpha}(s)$. To this end, let us analyze  what effect on the
constructions of Section \ref{section3} would be produced by the passing to the
general {\it tronqu\'ee} solution $u^{(tronq)}_{\alpha}(s)$, i.e., by lifting
the restriction $b_{1} = 1$. We already saw that the considerations
and results of Section \ref{section2} are not affected. The situation with
Section \ref{section3}, i.e., with the analysis of the Painlev\'e XXXIV
RH problem as $s \to -\infty$ is different.
In what follows, we shall analyze the {\it tronqu\'ee}-RH problem assuming that
\begin{equation}\label{bposit}
b >0.
\end{equation}

There are no changes in the basic three transformations of Section \ref{section3},
except that the factorization \eqref{factorization} now is
\begin{align} \label{factorization2}
      \begin{pmatrix} e^{-th(z)} & b
    \\ 0 &  e^{th(z)}\end{pmatrix} =
  \begin{pmatrix} 1 & 0
    \\ b^{-1} e^{th(z)} &  1\end{pmatrix}
   \begin{pmatrix} 0& b
    \\ -b^{-1} &  0\end{pmatrix}
 \begin{pmatrix} 1 & 0
    \\ b^{-1} e^{-th(z)} &  1\end{pmatrix},
\end{align}
and the transformation \eqref{Ydef} is modified accordingly.
It leads us to the following RH problem for the matrix function $C_{\alpha}(z)$.
\begin{enumerate}
\item[\rm (a)] $C_{\alpha} : \mathbb{C} \setminus \Sigma_C \to
    \mathbb C^{2\times 2}$ is analytic.
\item[\rm (b)] $C_{\alpha,+}(z) = C_{\alpha,-}(z)
    \begin{pmatrix} 1 & be^{-2t g(z)} \\ 0 & 1 \end{pmatrix}$,
    for $z \in (1,\infty)$,

    $C_{\alpha,+}(z) = C_{\alpha,-}(z)
    \begin{pmatrix} 0& b
    \\ -b^{-1} & 0\end{pmatrix}$,
    for $z \in (0,1)$ ,

    $C_{\alpha,+}(z) = C_{\alpha,-}(z)
    \begin{pmatrix} 0& 1
    \\ -1& 0\end{pmatrix}$,
    for $z \in \Sigma_3$,

    $C_{\alpha,+}(z) = C_{\alpha,-}(z)
    \begin{pmatrix} 1 & 0 \\  b^{-1}e^{2tg(z)} & 1 \end{pmatrix}$,
    for $z \in \Sigma_u \cup \Sigma_d$,

    $C_{\alpha,+}(z) = C_{\alpha,-}(z)
    \begin{pmatrix} 1 & 0 \\ e^{2\alpha \pi i + 2t g(z)} & 1 \end{pmatrix}$,
    for $z \in \Sigma_2$,

    $C_{\alpha,+}(z) = C_{\alpha,-}(z)
    \begin{pmatrix} 1 & 0 \\ e^{-2\alpha \pi i+2t g(z)} & 1 \end{pmatrix}$,
    for $z \in \Sigma_4$.
\item[\rm (c)]
$C_{\alpha}(z) =
\left(I + O\left(\frac{1}{z}\right)\right)
z^{-\sigma_3/4}
    \frac{1}{\sqrt{2}} \begin{pmatrix} 1 & i \\ i & 1 \end{pmatrix}$
    as $z \to \infty$.
\item[\rm (d)]
$C_{\alpha}(z) = O\begin{pmatrix} |z|^{\alpha} & |z|^{\alpha} \\
    |z|^{\alpha} & |z|^{\alpha}
\end{pmatrix}$ as $z \to 0$, if $-1/2 < \alpha < 0$; and

    $C_{\alpha}(z) = \left\{ \begin{array}{ll}
    O\begin{pmatrix} |z|^{\alpha} & |z|^{-\alpha} \\
    |z|^{\alpha} & |z|^{-\alpha} \end{pmatrix}
    & \text{as $z \to 0$ with $z \in (\Omega_1 \cup \Omega_4) \setminus (\Omega_u \cup \Omega_d)$}, \\[10pt]
    O\begin{pmatrix} |z|^{-\alpha} & |z|^{-\alpha} \\
    |z|^{-\alpha} & |z|^{-\alpha} \end{pmatrix}
    & \text{as $z \to 0$ with $z \in \Omega_2 \cup \Omega_3 \cup \Omega_u \cup \Omega_d$},
    \end{array} \right.$
    if $\alpha \geq 0$.
\end{enumerate}
The contour for this RH problem is the same as before, i.e., the one
depicted in Figure~\ref{figure5}.

We can at once make two important observations. Firstly,
 the neighborhood of the point $z =0$ will contribute
to the asymptotic analysis (as, in fact, it has in the
case $b =1$) and hence we should expect a change in the asymptotics
\eqref{uminusinfinity} and appearance in it of an explicit dependence
of the parameter $b$. Secondly, the inequality $b \neq 1$ yields
serious alterations in the constructions of the parametrices $P^{(\infty)}$
and $P_{\alpha}^{(0)}$.

The RH problem for the global parametrix $P^{(\infty)}$ now reads.
\begin{enumerate}
\item[\rm (a)] $P^{(\infty)} : \mathbb{C} \setminus (-\infty, 1] \to
    \mathbb C^{2\times 2}$ is analytic.
\item[\rm (b)] $P^{(\infty)}_{+}(z) = P^{(\infty)}_{-}(z)
    \begin{pmatrix} 0 & 1 \\ -1 & 0\end{pmatrix}$,
    for $z \in (-\infty,0)$,

$P^{(\infty)}_{+}(z) = P^{(\infty)}_{-}(z)
    \begin{pmatrix} 0 & b \\ -b^{-1} & 0\end{pmatrix}$,
    for $z \in (0,1)$,
\item[\rm (c)] $P^{(\infty)}(z) =
\left(I + O\left(\frac{1}{z}\right)\right)
z^{-\sigma_3/4}
    \frac{1}{\sqrt{2}} \begin{pmatrix} 1 & i \\ i & 1 \end{pmatrix}$
    as $z \to \infty$.
\end{enumerate}
This RH problem still can be solved explicitly. In fact, it
is now similar to the global parametrix from Section \ref{section2}.
The solution is given by the equation,
\begin{equation}\label{Pinftyb}
P^{(\infty)}(z)
= E (z-1)^{-\sigma_3/4} \frac{1}{\sqrt{2}} \begin{pmatrix} 1 & i \\ i & 1\end{pmatrix}
    \left(\frac{(z-1)^{1/2} + i}{(z-1)^{1/2} - i} \right)^{\beta \sigma_3}
\end{equation}
where
\[
\beta = \frac{i}{2\pi} \log b, \qquad E = \begin{pmatrix} 1 & 0 \\ 2\beta & 1\end{pmatrix},
\]
and the branches of the arguments are fixed by the inequalities,
\[
-\pi < \arg (z-1)< \pi, \quad 0 < \arg \bigg( \frac{(z - 1)^{1/2} + i }{(z - 1)^{1/2} - i} \bigg) < \pi.
\]

The  construction of the local parametrix $P^{(0)}_{\alpha}$ now involves,
instead of the Bessel model RH problem $\Phi_{\alpha}^{(Bes)}$,
the function $\Phi_{\alpha, \beta}^{(CHF)}$
which satisfies a RH problem with the same contour $\Gamma$ as $\Phi_{\alpha}^{(Bes)}$,
 and with the same jump matrices
in the left half-plane, see Figure \ref{figure7},
while in the right half-plane one has the new jump matrices,
\[
\begin{pmatrix} 1 & 0 \\ e^{2\pi i \beta} & 1\end{pmatrix},
\quad \begin{pmatrix} 0 & e^{-2\pi i \beta} \\ -e^{2\pi i \beta}& 0\end{pmatrix},
\quad\text{and}\quad \begin{pmatrix} 1 & 0 \\ e^{2\pi i \beta} & 1\end{pmatrix},
\]
on the rays $\arg \zeta = \pi/3$, $\arg \zeta = 0$, and $\arg \zeta = -\pi/3$, respectively.
The jump contour is  as shown in Figure~\ref{figure7}, but on straight rays extending to infinity.

\begin{figure}[th]
\centering
\unitlength 1pt
\linethickness{0.5pt}
\begin{picture}(200,160)(-100,-80)
   \put(0,0){\line(1,0){100}}
   \put(0,0){\line(-1,0){100}}
   \qbezier(0,0)(-25,40)(-50,80)
   \qbezier(0,0)(-25,-40)(-50,-80)
   \qbezier(0,0)(30,50)(50,80)
   \qbezier(0,0)(30,-50)(50,-80)
   \qbezier(0,-80)(0,0)(0,80)
   \qbezier(50,0)(50,0)(45,5)   \qbezier(50,0)(50,0)(45,-5)
   \qbezier(-60,0)(-60,0)(-65,5)   \qbezier(-60,0)(-60,0)(-65,-5)
   \qbezier(-25,40)(-25,40)(-23.5,47.5)   \qbezier(-25,40)(-25,40)(-33,42)
   \qbezier(-25,-40)(-25,-40)(-23.5,-47.5)   \qbezier(-25,-40)(-25,-40)(-33,-42)
   \qbezier(25,40)(25,40)(17,38)   \qbezier(25,40)(25,40)(26.5,32.5)
   \qbezier(25,-40)(25,-40)(17,-38)   \qbezier(25,-40)(25,-40)(26.5,-32.5)
   \qbezier(0,50)(1,49)(5,45)   \qbezier(0,50)(-1,49)(-5,45)
   \qbezier(0,-50)(1,-49)(5,-45)   \qbezier(0,-50)(-1,-49)(-5,-45)
   \put(-2,-12){$0$}
   \put(0,0){\circle*{4}}
   \put(100,3){$\begin{pmatrix} 0& e^{-\pi i\beta}  \\ -e^{\pi i\beta} & 0\end{pmatrix}$}
   \put(-120,60){$\begin{pmatrix} 1 & 0 \\  e^{2\pi i\alpha - \pi i\beta} & 1 \end{pmatrix}$}
   \put(-180,3){$\begin{pmatrix} 0& e^{\pi i\beta}  \\ -e^{-\pi i\beta} & 0\end{pmatrix}$}
   \put(-122,-60){$\begin{pmatrix} 1 & 0 \\ e^{-2\pi i\alpha  -\pi i\beta} & 1 \end{pmatrix}$}
   \put(47,65){$\begin{pmatrix} 1 & 0 \\  e^{-2\pi i\alpha+\pi i\beta} & 1 \end{pmatrix}$}
   \put(45,-69){$\begin{pmatrix} 1 & 0 \\ e^{2\pi i\alpha+\pi i\beta} & 1 \end{pmatrix}$}
   \put(-10,80){$e^{\pi i \alpha \sigma_3}$}
   \put(-10,-90){$e^{\pi i \alpha \sigma_3}$}
    \put(15,-70){I}
    \put(70,-30){II}
    \put(70,30){III}
    \put(15,70){IV}
    \put(-25,70){V}
    \put(-90,30){VI}
    \put(-90,-30){VII}
    \put(-30,-70){VIII}
   \end{picture}
   \caption{Contours and jump matrices for the RH problem for
   $\widetilde{\Phi}_{\alpha, \beta}^{(CHF)}$.}\label{figure11}
\end{figure}
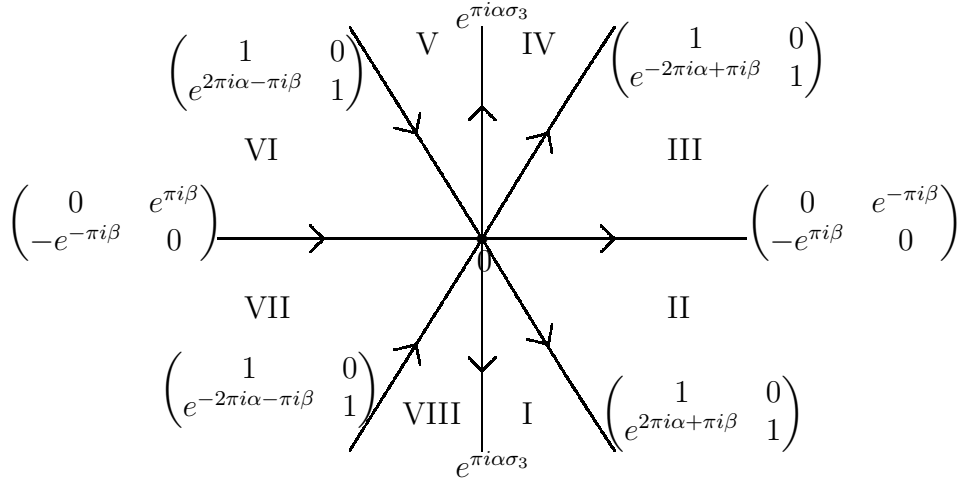

The function
$\Phi_{\alpha, \beta}^{(CHF)}$, in turn, admits the following representation,
\begin{equation}\label{CHFtilde}
    \Phi_{\alpha, \beta}^{(CHF)}(\zeta) =
    \widetilde{\Phi}_{\alpha, \beta}^{(CHF)}(\zeta)
    \left\{ \begin{array}{ll}
     e^{ \left(\frac{\pi i\beta}{2} +\pi i\alpha\right)\sigma_3},
        & \text{ if $\Re \zeta > 0$, $\Im \zeta > 0$}, \\
    e^{\frac{\pi i\beta}{2}\sigma_3} & \text{ if $\Re \zeta < 0$}, \\
     e^{\left(\frac{\pi i\beta}{2} - \pi i\alpha\right)\sigma_3},
        & \text{ if $\Re \zeta > 0$, $\Im \zeta < 0$},
    \end{array} \right.
\end{equation}
where the matrix function  $\widetilde{\Phi}_{\alpha, \beta}^{(CHF)}(\zeta)$
is the solution of the RH problem depicted in Figure~\ref{figure11}.

It is shown in \cite{DIK} that the RH problem for
$\widetilde{\Phi}_{\alpha, \beta}^{(CHF)}$ supplemented by the
proper representation at $\zeta =0$ (inherited  from \eqref{connection0P34}
and \eqref{connection1P34}) and the asymptotic condition,
\begin{equation}\label{CHFinfty}
\widetilde{\Phi}_{\alpha, \beta}^{(CHF)}(\zeta)
= \left(I + O\left(\frac{1}{\zeta}\right)\right)\zeta^{-\beta\sigma_{3}}
    e^{-i \zeta \sigma_3}\qquad \text{as}\quad \zeta \to \infty,
\quad 0 < \arg\zeta < \frac{\pi}{2},
\end{equation}
is uniquely solvable; moreover, it admits an explicit solution in terms of
the confluent hypergeometric functions $\psi(a,c;\zeta)$ with the parameters,
\[
a = \alpha + \beta,\qquad c = 1+2\alpha.
\]
Indeed, the solution $\widetilde{\Phi}_{\alpha, \beta}^{(CHF)}$ is described
by the following formulae.

Define on the complex plane, cut along the negative imaginary axis,
the matrix function,
\begin{multline}\label{PsiConfl1}
\Psi_{0}(\zeta):=\left(\begin{matrix}
2^{\alpha}\zeta^{\alpha}\psi(\alpha+\beta,1+2\alpha, 2e^{\frac{i\pi}{2}}
\zeta)e^{i\pi\left(2\beta+\frac{3\alpha}{2}\right)}e^{-i\zeta} \cr
-2^{-\alpha}\zeta^{-\alpha}
\psi(1-\alpha+\beta,1-2\alpha, 2e^{\frac{i\pi}{2}}\zeta)e^{i\pi\left(\beta-\frac{7\alpha}{2}\right)}e^{-i\zeta}
{\Gamma(1+\alpha+\beta)\over\Gamma(\alpha-\beta)}
\end{matrix}
\right.\\
\left.
\begin{matrix}
-2^{\alpha}\zeta^{\alpha}
\psi(1+\alpha-\beta,1+2\alpha, 2e^{-\frac{i\pi}{2}}\zeta)e^{i\pi\left(\beta+\frac{3\alpha}{2}\right)}e^{i\zeta}
{\Gamma(1+\alpha-\beta)\over\Gamma(\alpha+\beta)}
\cr
2^{-\alpha}\zeta^{-\alpha}
\psi(-\alpha-\beta,1-2\alpha, 2^{-\frac{i\pi}{2}}\zeta)e^{-\frac{3i\pi}{2}\alpha}
e^{i\zeta}\end{matrix}\right),
\end{multline}
where the branches of the multi-valued functions in the right hand side
of the equation (including the confluent hypergeometric function) are
fixed by the condition,
\[
-\frac{\pi}{2} < \arg \zeta < \frac{3\pi}{2}.
\]
We use I--VIII to denote the eight sectors as in Figure \ref{figure11}. The function
$\widetilde{\Phi}_{\alpha, \beta}^{(CHF)}(\zeta)$ is given then by the equations,
\begin{multline}
 \widetilde{\Phi}_{\alpha, \beta}^{(CHF)}(\zeta) = 2^{\beta\sigma_{3}}e^{\frac{i\pi}{2}\beta\sigma_{3}}
 \begin{pmatrix} e^{-i\pi(\alpha +2\beta)}&0  \\ 0& e^{i\pi(2\alpha +\beta)}\end{pmatrix}
 \\ \label{CHFtilde2}
 \times \Psi_{0}(\zeta) \left\{ \begin{array}{ll}
     \begin{pmatrix} 0&- e^{-\pi i\beta}  \\ e^{\pi i\beta} & 0\end{pmatrix}
        & \text{ if $\zeta \in I$}, \\[10pt]
     \begin{pmatrix} e^{2\pi i \alpha}&- e^{-\pi i\beta}  \\ e^{\pi i\beta} & 0\end{pmatrix}
        & \text{ if $\zeta \in II$}, \\[10pt]
      I & \text{ if $\zeta \in III$}, \\
     \begin{pmatrix}1&0  \\ e^{\pi i(\beta-2\alpha)} & 1\end{pmatrix}
        & \text{ if $\zeta \in IV$}, \\[10pt]
     \begin{pmatrix} e^{\pi i \alpha}&0 \\ e^{\pi i(\beta-\alpha)} & e^{-\pi i \alpha}\end{pmatrix}
        & \text{ if $\zeta \in V$}, \\[10pt]
      \begin{pmatrix} e^{\pi i \alpha}&0 \\ 2i\sin \pi(\beta - \alpha)& e^{-\pi i \alpha}\end{pmatrix}
        & \text{ if $\zeta \in VI$},\\[10pt]
       \begin{pmatrix}0&-e^{i\pi(\alpha + \beta)} \\ e^{-\pi i(\alpha+\beta)} &
        -2ie^{i\pi \beta}\sin \pi(\beta - \alpha)\end{pmatrix}
        & \text{ if $\zeta \in VII$},\\[10pt]
        \begin{pmatrix}e^{-i\pi\alpha}&-e^{i\pi(\alpha + \beta)} \\ e^{\pi i(\beta-3\alpha)} &
        -2ie^{i\pi \beta}\sin \pi(\beta - \alpha)\end{pmatrix}
        & \text{ if $\zeta \in VIII$}.
          \end{array} \right.
\end{multline}
It is a straightforward though a bit involved calculation to check that the function $\widetilde{\Phi}_{\alpha, \beta}^{(CHF)}$
defined by (\ref{PsiConfl1})-(\ref{CHFtilde2}) does indeed satisfy the jump conditions indicated in Figure \ref{figure11}.
By a  direct calculation, one can also establish the following asymptotic  behavior of
$\widetilde{\Phi}_{\alpha, \beta}^{(CHF)}(\zeta)$ as $\zeta \to \infty$.
\begin{multline}
 \widetilde{\Phi}_{\alpha, \beta}^{(CHF)}(\zeta) = \left(I + \frac{m_{\tilde{\Phi}}}{\zeta} + O\left(\frac{1}{\zeta^2}\right)\right)
\zeta^{-\beta\sigma_{3}}e^{-i\zeta\sigma_{3}}
 \\ \label{CHFtilde2as}
 \times  \left\{ \begin{array}{ll}
     I & \text{if } 0 < \arg\zeta < \frac{\pi}{2}, \\[5pt]
     e^{i\pi\alpha\sigma_{3}}& \text{if } \frac{\pi}{2} < \arg\zeta < \pi, \\[5pt]
      \begin{pmatrix}0&-e^{i\pi(\alpha + \beta)}  \\  e^{-i\pi(\alpha + \beta)} &0\end{pmatrix}
        & \text{if } \pi < \arg\zeta < \frac{3\pi}{2}, \\[10pt]
        \begin{pmatrix}0&-e^{-i\pi \beta}  \\  e^{i\pi \beta}&0\end{pmatrix}
        & \text{if } -\frac{\pi}{2} < \arg\zeta < 0,
    \end{array} \right.
\end{multline}
where
\begin{equation}\label{m1CHF}
m_{\tilde{\Phi}} = \begin{pmatrix}\frac{i}{2}(\beta^2-\alpha^2)&
-\frac{i}{2}\frac{\Gamma(1+\alpha-\beta)}{\Gamma(\alpha+\beta)}e^{4\pi i \alpha + i\pi\beta}  \\[10pt]
\frac{i}{2}\frac{\Gamma(1+\alpha+\beta)}{\Gamma(\alpha-\beta)}e^{-4\pi i \alpha - i\pi\beta} &
\frac{i}{2}(\alpha^2-\beta^2)\end{pmatrix}.
\end{equation}

The asymptotic formulae (\ref{CHFtilde2as}) indicate that the local parametrix $P_{\alpha,\beta}^{(0)}(z)$
takes again the form,
\begin{equation}
P_{\alpha,\beta}^{(0)}(z) = E_{\alpha,\beta}(z)\Phi_{\alpha, \beta}^{(CHF)}(tf(z))e^{tg(z)\sigma_{3}},
\end{equation}
where the change-of-variable function $f(z)$ is exactly the same as before, i.e., as in the case $b=1$,
while the holomorphic at $z =0$ matrix valued function $E_{\alpha,\beta}(z)$ is defined by
the formula which is slightly more complicated than the previous equations (\ref{defEalpha}).
Indeed, this time we have,
\begin{align} \label{defEalpha1}
    E_{\alpha,\beta}(z) = \left\{ \begin{array}{ll}
    P^{(\infty)}(z)(tf(z))^{\beta\sigma_{3}} e^{2it\sigma_3/3}e^{-i\pi\left(\alpha +\frac{\beta}{2}\right)\sigma_{3}},
    & \quad \Im z > 0, \\[10pt]
    P^{(\infty)}(z) \begin{pmatrix} 0 & 1 \\ -1 & 0 \end{pmatrix}
    (tf(z))^{\beta\sigma_{3}} e^{2it\sigma_3/3}e^{-i\pi\left(\alpha +\frac{\beta}{2}\right)\sigma_{3}}, & \quad \Im z < 0,
 \end{array} \right.
        \end{align}
with $0 <\arg f(z) < 2\pi$, and, in particular,
\begin{equation}\label{Eab0}
E_{\alpha,\beta}(0) = \frac{1}{\sqrt{2}}\begin{pmatrix}1&1\\ 2\beta -1 & 2\beta +1\end{pmatrix}
(4t)^{\beta\sigma_{3}}e^{\frac{2it}{3} -i\pi\left(\frac{1}{4} +\alpha - \frac{\beta}{2}\right)\sigma_{3}}.
\end{equation}

Having constructed the local parametrix at $z =0$, the further arguments are identical
to the ones we have used in the case $b=1$. As a result, we arrive at the following  representation for
the solution $C_{\alpha}(z)$  to  the ``master'' $C$ -RH problem (cf. (\ref{fromCtoD})).
\begin{align} \nonumber
    C_{\alpha}(z) & = D_{\alpha, \beta}(z) P_{\alpha, \beta}^{(0)}(z) \\
        & = D_{\alpha, \beta}(z) E_{\alpha, \beta}(z)
        \Phi_{\alpha, \beta}^{(CHF)}(t f(z))
                e^{t g(z) \sigma_3}, \label{fromCtoD2}
                \end{align}
where $D_{\alpha, \beta}(z) = I + O((-s)^{-3/2})$ as $ s \to -\infty$. Hence,
similar to the case $b =1$, we have from (\ref{usolution4C}) and (\ref{fromCtoD2})
that
\begin{align} \nonumber
   u(s) &=  \frac{i}{\sqrt{-s}}\left[E_{\alpha,\beta}(0)
        \lim_{z \to 0}
        \left( z \left(\frac{d}{dz} \Phi_{\alpha, \beta}^{(CHF)}(t f(z)) \right)
         \left( \Phi_{\alpha, \beta}^{(CHF)}(t f(z))
        \right)^{-1}\right)E^{-1}_{\alpha,\beta}(0) \right]_{12} + O(s^{-2}) \\
    & = \nonumber
      \frac{i}{\sqrt{-s}}\left[E_{\alpha,\beta}(0)
        \lim_{\zeta \to 0}
        \left( \zeta \left(\frac{d}{d\zeta} \Phi_{\alpha, \beta}^{(CHF)}(\zeta) \right)
         \left( \Phi_{\alpha, \beta}^{(CHF)}(\zeta)
        \right)^{-1}\right)E^{-1}_{\alpha,\beta}(0) \right]_{12} + O(s^{-2}) \\
    & = \label{ubasmin1}
   \frac{i}{\sqrt{-s}}\left[E_{\alpha,\beta}(0)
        \lim_{\zeta \to 0}
        \left( \zeta \left(\frac{d}{d\zeta}\widetilde{\Phi}_{\alpha, \beta}^{(CHF)}(\zeta) \right)
         \left( \widetilde{\Phi}_{\alpha, \beta}^{(CHF)}(\zeta)
        \right)^{-1}\right)E^{-1}_{\alpha,\beta}(0) \right]_{12} + O(s^{-2}),
\end{align}
as $s \to -\infty$. The last equality follows from the fact that the  matrix function
$\widetilde{\Phi}_{\alpha, \beta}^{(CHF)}(\zeta)$ differs from
the matrix function $\Phi_{\alpha, \beta}^{(CHF)}(\zeta)$ only by a non-singular piecewise constant
right matrix multiplier, see (\ref{CHFtilde}).

The function $\widetilde{\Phi}_{\alpha, \beta}^{(CHF)}(\zeta)$ solves the Riemann-Hilbert
problem whose jump conditions are  depicted in Figure \ref{figure11},
the asymptotic behavior as $\zeta \to \infty$ is indicated in (\ref{CHFinfty}), and the
branching singularity at $\zeta = 0$ is of type described in (\ref{connection0P34})
and  (\ref{connection1P34}). Indeed,  if $\alpha - \frac{1}{2} \not\in \mathbb N_0$, then
\begin{equation}\label{connection0CHF}
   \widetilde{\Phi}_{\alpha, \beta}^{(CHF)}(\zeta) = B(\zeta)
    \begin{pmatrix} \zeta^{\alpha} & 0 \\ 0 & \zeta^{-\alpha} \end{pmatrix}Q,
\end{equation}
where $B$ is analytic at $\zeta =0$.
If $\alpha \in \frac{1}{2} + \mathbb N_0$, then there exists a constant $\kappa$ such that
\begin{equation} \label{connection1CHF}
    \widetilde{\Phi}_{\alpha, \beta}^{(CHF)}(\zeta) = B(\zeta)
    \begin{pmatrix} \zeta^{\alpha} & \kappa \zeta^{\alpha}\log \zeta \\ 0 & \zeta^{-\alpha} \end{pmatrix}Q,
   \end{equation}
where $B$ is again analytic at $\zeta =0$. The right matrix
multipliers $Q$ in (\ref{connection0CHF}) and
(\ref{connection1CHF}) are piecewise constant matrix functions; in
fact, they are  constant in the eight sectors I--VIII. The
matrices $Q$ can be written down explicitly using either  the
general algebraic properties of the Riemann-Hilbert problem, or
the explicit formulae (\ref{PsiConfl1}) - (\ref{CHFtilde2}) for
the function $\widetilde{\Phi}_{\alpha, \beta}^{(CHF)}(\zeta)$.
However, we won't do this. The important feature of the matrices
$Q$ as well as of all the jump matrices of the
$\widetilde{\Phi}_{\alpha, \beta}^{(CHF)}$ - RH problem, is that
they are constant with respect to $\zeta$. Therefore, we can
exploit the standard arguments of the theory of integrable systems
(see e.g. \cite{JMU}, \cite{FT}; see also \cite{FIKN},   Chapters
2, 3) and conclude that
\begin{equation}\label{logder1}
       \left(\frac{d}{d\zeta}\widetilde{\Phi}_{\alpha, \beta}^{(CHF)}(\zeta)\right)
         \left( \widetilde{\Phi}_{\alpha, \beta}^{(CHF)}(\zeta)
        \right)^{-1}
\equiv -i\zeta\sigma_{3}  +\frac{A_{-1}}{\zeta},
\end{equation}
where
\begin{equation}\label{A1}
A_{-1} = -\beta\sigma_{3} + i[\sigma_{3}, m_{\tilde{\Phi}}]
= \begin{pmatrix}-\beta&
\frac{\Gamma(1+\alpha-\beta)}{\Gamma(\alpha+\beta)}e^{4\pi i \alpha + i\pi\beta}  \\[5pt]
\frac{\Gamma(1+\alpha+\beta)}{\Gamma(\alpha-\beta)}e^{-4\pi i \alpha - i\pi\beta} &
\beta\end{pmatrix}.
\end{equation}
Indeed, the $\zeta$-independence of all  the jump matrices
implies that the function{\footnote{All the jump matrices of
the $\widetilde{\Phi}_{\alpha, \beta}^{(CHF)}$ - RH problem are having the unit determinant.
Together with the asymptotic condition (\ref{CHFinfty}) this yield
the identity, $\det\widetilde{\Phi}_{\alpha, \beta}^{(CHF)}(z) \equiv 1$ and
therefore the holomorphic invertability of the matrix $\widetilde{\Phi}_{\alpha, \beta}^{(CHF)}(\zeta)$
for all $\zeta$. }},
\[
 A(\zeta):= \left( \frac{d}{d\zeta} \widetilde{\Phi}_{\alpha, \beta}^{(CHF)}(\zeta) \right)
         \left( \widetilde{\Phi}_{\alpha, \beta}^{(CHF)}(\zeta)
        \right)^{-1},
\]
is analytic on ${\Bbb C}\setminus \{0\}$. Moreover, from the behavior of
$\widetilde{\Phi}_{\alpha, \beta}^{(CHF)}(\zeta)$ at $\zeta=0$ (see (\ref{connection0CHF}),
(\ref{connection1CHF})) and at $\zeta=\infty$ (see (\ref{CHFtilde2as})),  we
conclude that $A(\zeta)$ is in fact a rational function which has the only simple pole
at $\zeta=0$ and such that $A(\infty) = i\sigma_{3}$. Hence, by Liouville's theorem,
\[
A(\zeta) = -i\sigma_{3} + \frac{A_{-1}}{\zeta}.
\]
This proves (\ref{logder1}). Equation (\ref{A1}) is obtained by
substituting  the asymptotic expansion (\ref{CHFtilde2as})
into the left hand side of  (\ref{logder1}) and
equating the terms of the order
 {\footnote{Of course, equation (\ref{logder1})
can be deduced via the direct differentiation of (\ref{PsiConfl1})
and use of the classical confluent hypergeometric equation. This
derivation though is much more involved then the one we just used.}}
$O(\zeta^{-1})$.

Relation (\ref{logder1}) allows us to re-write the asymptotic formula (\ref{ubasmin1})
as
\begin{equation}\label{ubasmin2}
    u(s) = \frac{i}{\sqrt{-s}}\Bigl[E_{\alpha,\beta}(0)A_{-1}E^{-1}_{\alpha,\beta}(0)\Bigr]_{12}
    + O(s^{-2}),\quad s \to -\infty.
\end{equation}
Plugging in here (\ref{Eab0}) and (\ref{A1}) yields the asymptotic
equation,
\begin{equation}\label{ubasmin3}
    u(s) = \frac{i}{\sqrt{-s}}
    \left[\beta - \frac{i}{2}\frac{\Gamma(1+\alpha -\beta)}{\Gamma(\alpha + \beta)}e^{i\theta(s)}
    - \frac{i}{2}\frac{\Gamma(1+\alpha +\beta)}{\Gamma(\alpha - \beta)}e^{-i\theta(s)}\right]
    + O\left(\frac{1}{s^2}\right),\quad s \to -\infty,
\end{equation}
where
\[
\theta(s) = \frac{4t}{3} - 2i\beta\log{8t} - \alpha \pi, \qquad t = (-s)^{3/2}.
\]
Taking into account that we have assumed that $b>0$ and hence that $\beta$ is
pure imaginary, we transform (\ref{ubasmin3}) into our final asymptotic
representation  for the (real-valued) {\it tronqu\'ee} solutions of the
thirty fourth Painlev\'e equation as $s \to -\infty$,
\begin{equation}\label{tronqminus}
u^{(tronq)}(s) = \frac{\beta_{0}}{\sqrt{-s}} +\frac{\sqrt{\alpha^2
+ \beta^{2}_{0}}}{\sqrt{-s}} \cos \left(\frac{4}{3}(-s)^{3/2}
-3\beta_{0}\log(-s) + \chi\right) + O(s^{-2}),
\end{equation}
where
\[
    \beta_{0} \equiv i\beta = -\frac{1}{2\pi}\log b,
\]
and the phase $\chi$ is given explicitly in terms of
$\beta_{0}$ and $\alpha$,
\begin{equation}\label{phaseb}
\chi = -\alpha \pi -6\beta_{0}\log 2 + \arg (\alpha + i \beta_{0})
+ 2\arg \Gamma(\alpha + i \beta_{0}).
\end{equation}

\begin{remark} The above asymptotic analysis of the Painlev\'e XXXIV Riemann-Hilbert
problem does not need the condition $b>0$ and hence the pure
imaginary $\beta$. The necessary restriction is in fact,
\[
|\Re \beta | < \frac{1}{2}.
\]
With this restriction, we again arrive to the more general formula (\ref{ubasmin3}),
but with the change of the error term
\[
 O\left(\frac{1}{s^2}\right)
\]
to the term
\[
 O\left(\frac{1}{s^{2-3|\Re\beta|}}\right).
\]
\end{remark}

The asymptotic analysis of the general {\it tronqu\'ee}
solution which we just performed makes a strong case in favor of the fact  that
the solution $u_{\alpha}(s)$ is characterized uniquely
by the asymptotic conditions of Theorem \ref{theorem2}.
Indeed, using the transformation formula \eqref{p34top2},
we conclude that every {\it tronqu\'ee} solution of
Painlev\'e XXXIV, i.e., one behaving at $s \to +\infty$
as \eqref{uplusinfinity}, maps to a {\it tronqu\'ee} solution of
Painlev\'e II.  
This means that the Painlev\'e II
Stokes multipliers  must be as in \eqref{specP2a2free}, because
otherwise the asymptotics is oscillatory for $s \to -\infty$ and
elliptic in the sectors $2\pi/3 < \arg s < \pi $ and $\pi < \arg s < 4\pi/3 $
(see \cite{Kap1} and \cite{IKap1}). This in turns means
that the asymptotic condition \eqref{uplusinfinity} selects
the {\it tronqu\'ee} family \eqref{btrongp34} of solutions.
As we have already shown, for the positive $b$, in fact for
all $b$ such that $|\arg b| < \pi$ it must be
exactly $b =1$ in order for the solution to behave at $-\infty$ as
it is indicated in \eqref{uminusinfinity}. It is natural to expect that
in the case $b < 0$ the ``minus infinity''  asymptotics of the solution
is different from \eqref{uminusinfinity}.

As it was mentioned in the introduction, we actually conjecture that the
asymptotics  \eqref{uminusinfinity} alone fixes the solution uniquely.
The supporting arguments are based on the fact that the substitution
of the asymptotic ansatz \eqref{tronqminus} into the formula \eqref{p34top2}
yields the asymptotic representation for the Painlev\'e II function,
\begin{equation}\label{p2gen}
q(s)\sim \sqrt{\frac{s}{2}}\cot
 \left(\frac{\sqrt{2}}{3}s^{3/2} -\frac{3}{2}\beta_{0}\log s + \hat{\chi}\right)
 \qquad \text{as } s \to \infty.
\end{equation}
This asymptotics is consistent with the general results of \cite{Kap1}
for real-valued solutions of Painlev\'e II which provide
the explicit formulae relating the Painlev\'e II monodromy data
$a_{1}$, $a_{2}$ and $a_{3}$ and the asymptotic parameters $\beta_{0}$
and $\hat{\chi}$. This, in turn, allows us to establish
a one-to-one correspondence between the Painlev\'e XXXIV monodromy data
$b_{1}$, $b_{2}$ and $b_{4}$ and the asymptotic parameters $\beta_{0}$ and
 $\chi$. Hence the uniqueness of the real-valued solution of
 the Painlev\'e XXXIV equation \eqref{painleve34} with the asymptotic
 condition \eqref{uminusinfinity}.


\begin{thebibliography}{99}
\bibitem{AS} M. Abramowitz and I. Stegun,
    Handbook of Mathematical Functions,
    Dover Publications, New York, 1992. Reprint of the 1972 edition.

\bibitem{BBIK}
    P. Bleher, A. Bolibruch, A. Its, and A. Kapaev,
    Linearization of the P34 equation of Painlev\'e-Gambier,
    unpublished manuscript.

\bibitem{CKV}
    T. Claeys, A.B.J. Kuijlaars and M. Vanlessen,
    Multi-critical unitary random matrix ensembles and the general
    Painlev\'e II equation,
    Ann. Math. 167 (2008), 601--642.

\bibitem{Deift}
    P. Deift,
    Orthogonal Polynomials and Random Matrices: A  Riemann-Hilbert Approach,
    Courant Lecture Notes 3, New York University, 1999.

\bibitem{DIK}
    P. Deift, A. Its, and I. Krasovsky,
    Toeplitz and Hankel determinants with Fisher-Hartwig singularities,
    in preparation.

\bibitem{DKMVZ2}
    P. Deift, T. Kriecherbauer, K.T-R McLaughlin, S. Venakides, and X. Zhou,
    Uniform asymptotics for polynomials orthogonal with respect to
    varying exponential weights and applications to universality
    questions in random matrix theory,
    Comm. Pure Appl. Math. 52 (1999), 1335--1425.

\bibitem{DZ}
    P. Deift and X. Zhou,
    A steepest descent method for oscillatory Riemann-Hilbert problems.
    Asymptotics for the MKdV equation,
    Ann. Math. 137 (1993), 295--368.


\bibitem{FT}
L.D. Faddeev and L.A. Takhtajan,
 {\it Hamiltonian Methods  in the Theory of Solitons},
  Springer Verlag, Berlin, Heidelberg, 1987.

\bibitem{FN}
    H. Flaschka and A.C. Newell,
    Monodromy and spectrum-preserving deformations I,
    Comm. Math. Phys. 76 (1980), 65--116.

\bibitem{FIKN}
    A.S. Fokas, A.R. Its, A.A. Kapaev, and V.Yu. Novokshenov,
    Painlev\'e Transcendents: The Riemann-Hilbert Approach,
    AMS  Mathematical Surveys and Monographs, vol. 128,
    Amer. Math. Society, Providence R.I., 2006.

\bibitem{Ince}
    E.L. Ince,
    Ordinary Differential Equations,
    Dover, New York, 1944.

\bibitem{IKap1}
    A.R. Its and A.A. Kapaev,
    The irreducibility of the second Painlev\'e equation and the isomonodromy method.
    In: Toward the exact WKB analysis of differential equations, linear or non-linear,
    C.J. Howls, T. Kawai, and Y. Takei, eds.,
    Kyoto Univ. Press, 2000, pp. 209--222.

\bibitem {IKap2}
    A.R. Its and A.A. Kapaev,
    Quasi-linear Stokes phenomenon for the second Painlev\'e transcendent,
    Nonlinearity 16 (2003), 363--386.

\bibitem{IK}
    A. Its and I. Krasovsky,
    Hankel determinant and orthogonal polynomials for the Gaussian weight with a jump,
    in ``Integrable Systems and Random Matrices'' (J. Baik et al., eds.),
    Contemporary  Mathematics 458, Amer. Math. Soc., Providence R.I. 2008,
    pp. 215--248.

\bibitem{IKO1}
    A.R. Its, A.B.J. Kuijlaars, and J. \"Ostensson,
    Critical edge behavior in unitary random matrix ensembles and
    the thirty fourth Painlev\'e transcendent,
    Internat. Math. Research Notices,  2008 Volume 2008: article ID rnn017, 67 pages.

\bibitem{JMU}
    M. Jimbo, T. Miwa, and K. Ueno,
    Monodromy preserving deformation of linear ordinary differential
    equations with rational coefficients,
    Physica D  2 (1980), 306--352.

\bibitem{Kap1}
    A.A. Kapaev,
    Global asymptotics of the second Painlev\'e transcendent,
    Phys. Lett. A, 167 (1992) 356--362.

\bibitem{Kap2}
    A.A. Kapaev,
    Quasi-linear Stokes phenomenon for the Hastings-McLeod
    solution of the second Painlev\'e equation,
    arXiv: nlin.SI/0410009.

\bibitem{KH}
    A.A. Kapaev and E. Hubert,
    A note on the Lax pairs for Painlev\'e equations,
    J. Phys. A  32 (1999), 8145--8156.

\bibitem{KV}
    A.B.J. Kuijlaars and M. Vanlessen,
    Universality for eigenvalue correlations at the origin of the spectrum,
    Comm. Math. Phys.  243  (2003), 163--191.

\bibitem{SS} N. Seiberg and D. Shih,
    Flux vacua and branes of the minimal superstring,
    J. High Energy Phys. 2005, no. 1, 055, 38 pp.

\bibitem{TW1}
    C.A. Tracy and H. Widom,
    Level-spacing distributions and the Airy kernel,
    Comm. Math. Phys. 159 (1994), 151--174.

\bibitem{TW2}
    C.A. Tracy and H. Widom,
    Airy kernel and Painlev\'e II,
    in ``Isomonodromic Deformations and Applications in Physics'',
    (J. Harnad and A. Its, eds),
    CRM Proc. Lecture Notes, 31, Amer. Math. Soc., Providence, RI, 2002.
    pp.~85--96.

\bibitem{Vanl}
    M. Vanlessen,
    Strong asymptotics of the recurrence coefficients of orthogonal
    polynomials associated to the generalized Jacobi weight,
    J. Approx. Theory  125  (2003), 198--237.
\end{thebibliography}
\end{document}